\documentclass{amsart}

\newcommand{\baton}[1]{\mathbb #1}
\newcommand{\C}{{\baton C}}
\newcommand{\E}{{\baton E}}
\newcommand{\N}{{\baton N}}
\newcommand{\Q}{{\baton Q}}
\newcommand{\R}{{\baton R}}
\newcommand{\T}{{\baton T}}
\newcommand{\Z}{{\baton Z}}
\DeclareMathOperator{\circlegroup}{S}
\renewcommand{\S}{\circlegroup}

\newcommand{\CA}{{\mathcal A}}
\newcommand{\CB}{{\mathcal B}}
\newcommand{\CC}{{\mathcal C}}

\newcommand{\CF}{{\mathcal F}}

\newcommand{\CN}{{\mathcal N}}

\newcommand{\CP}{{\mathcal P}}
\newcommand{\CM}{{\mathcal M}}
\newcommand{\CS}{{\mathcal S}}
\newcommand{\CX}{{\mathcal X}}
\newcommand{\CZ}{{\mathcal Z}}

\newcommand{\norm}[1]{\lVert #1\rVert}
\DeclareMathOperator{\id}{Id}
\newcommand{\one}{{\boldsymbol 1}}
\newcommand{\inv}{^{-1}}
\newcommand{\wh}{\widehat}
\DeclareMathOperator{\AP}{{\CA\CP}}
\DeclareMathOperator{\Nil}{\CN^2}
\DeclareMathOperator{\Isom}{Isom}
\DeclareMathOperator{\Av}{Av}
\DeclareMathOperator{\Sp}{Sp}

\newtheorem{theorem}{Theorem}
\newtheorem{lemma}{Lemma}
\newtheorem{proposition}{Proposition}
\newtheorem*{proposition*}{Proposition}
\newtheorem{corollary}{Corollary}

\theoremstyle{definition}
\newtheorem{definition}{Definition}
\newtheorem*{notation}{Notation}
\newtheorem*{convention}{Convention}

\theoremstyle{remark}
\newtheorem*{remark*}{Remark}
\newtheorem{claim}{Claim}
\theoremstyle{plain}

\newcommand{\uq}{{\bf q}}
\newcommand{\ua}{{\bf a}}
\newcommand{\ub}{{\bf b}}
\newcommand{\uc}{{\bf c}}

\newcommand{\ue}{{\bf e}}

\newcommand{\uomega}{\boldsymbol{\omega}}

\DeclareMathOperator{\BZ}{B(\Z)}

\renewcommand{\theenumi}{\roman{enumi}}

\begin{document}

\title{Analysis of two step nilsequences}
\author{Bernard Host and Bryna Kra}

\address{ Universit\'e Paris-Est, Laboratoire d'analyse et de
math\'ematiques appliqu\'ees, UMR CNRS 8050, 5 bd Descartes, 77454 Marne la
Vall\'ee Cedex 2, France}
\email{bernard.host@univ-mlv.fr}
\address{Department of Mathematics, Northwestern University,
2033 Sheridan Road,  Evanston, IL 60208-2730, USA}
\email{kra@math.northwestern.edu}
\thanks{The second author  was partially supported by
NSF grant DMS-0555250.}

\begin{abstract}
Nilsequences arose in the study 
of the multiple ergodic averages associated 
to Furstenberg's proof of Szemer\'edi's Theorem 
and have since played a role in problems in 
additive combinatorics.  Nilsequences are a generalization 
of almost periodic sequences and we study which 
portions of the classical theory for almost periodic 
sequences can be generalized for two step 
nilsequences.  We state and prove basic properties 
for $2$-step nilsequences and give a classification scheme 
for them.  
\end{abstract}
\maketitle

\section{Introduction}
Traditional Fourier analysis has been used with great success to study 
problems that are linear in nature.  For example, Roth used 
Fourier techniques (the circle method) 
to show that a set of integers with positive 
upper density contains infinitely many arithmetic progressions 
of length $3$ and van der Corput used such methods
to show that the primes contain infinitely many arithmetic progressions 
of length $3$.   
These questions can be reformulated into solving a certain 
linear equation in a subset of the integers, making them tractable 
to Fourier analysis; the linear nature of the equation
is a key ingredient in the use of 
Fourier techniques.  

Finer analysis is needed to study 
longer arithmetic progressions or more complicated patterns.  
This analysis is necessarily quadratic (or of higher order) 
in nature.  
The basic objects in such an analysis are the {\em nilsequences} 
(the precise definition is given in~\ref{def:nilsequ}). 
Nilsequences were introduced in~\cite{BHK} 
for studying some problems in ergodic theory 
and additive combinatorics that arose from the multiple ergodic 
averages introduced by Furstenberg~\cite{F} in his 
ergodic theoretic proof of Szemer\'edi's Theorem.  
Two-step nilsequences play a role in the recent work of 
Green and Tao (see~\cite{GT} and~\cite{GT2})  on 
arithmetic progressions of length $4$ 
in the primes; higher order nilsequences are conjectured 
to play an analogous role for longer arithmetic progressions.  

In spite of the recent uses of nilsequences in a variety of contexts, 
as of yet there is no systematic study of them.  
Even some of the basic properties are neither stated  nor proved.
In this paper, we carry out such an analysis for $2$-step 
nilsequences and hope that this clarifies their role and properties. 
In a forthcoming paper~\cite{HK6}, 
we give applications of the results in this 
paper.

Almost periodic sequences and their classically understood properties 
are a model we would like to emulate for nilsequences.  
Almost periodic sequences (which are exactly 
$1$-step nilsequences) arise naturally in understanding 
the correlation sequence 
$$(\int f\cdot f\circ T^n\,d\mu\colon n\in\Z) \ , 
$$
where $(X, \CX, \mu, T)$ is a measure preserving 
probability system and $f\in L^\infty(\mu)$.  
In a similar way, $2$-step 
nilsequences arise in~\cite{BHK} for understanding the double 
correlation sequence 
$$
(\int f\cdot f\circ T^n\cdot f\circ T^{2n}\,d\mu\colon n\in\Z) \ . 
$$
Unfortunately, there are fundamental differences 
between the $1$-step (almost periodic) and $2$-step 
nilsequences that make it impossible to carry over 
the classical results in a straightforward manner.  
Namely, we no longer have functions such as the 
trigonometric functions that form a basis.  

On the other hand, there is a class of nilsequences 
called {\em elementary nilsequences} (see Definition~\ref{def:elementary}) 
that are the building blocks for all nilsequences.  
Using these, we obtain a density result (Theorem~\ref{th:span2})
analogous to the density of trigonometric polynomials in $L^2$. 
Restricting to the elementary nilsequences, we can 
classify them such that each elementary nilsequence belongs 
to a class with a particularly simple representative of this class. 

Throughout this article, we restrict ourselves to $2$-step nilsequences. 
Some of the analysis here can be carried out for higher levels, 
but most of the results are relative to the previous level of almost 
periodic sequences.  The detailed analysis of higher 
order nilsequences begins with a firm understanding of the 
$2$-step ones. 

\section{Definitions and examples}

\subsection{Almost periodic sequences}

\begin{notation}
In general, we denote a sequence by 
$\ua=(a_n\colon n\in\Z)$. Given $\ua$, we use $\sigma\ua$ to 
denote the {\em shifted} sequence.  Thus for $k\in\Z$, by $\sigma^k\ua$, 
we mean the shifted sequence $\sigma^k\ua = 
(a_{n+k}\colon n\in\Z)$.

For $t\in\T:=\R/\Z$, we use the standard notation
$e(t)=\exp(2\pi i t)$. 
The {\em exponential linear sequence}  $\ue(t)$ 
is defined by $\ue(t)=(e(nt)\colon n\in\Z)$. 
\end{notation}

\begin{proposition}[and definition] 
\label{prop:AP}
For a bounded sequence 
$\ua=(a_n\colon n\in\Z)$, the following properties are equivalent:
\begin{enumerate}
\item
\label{item:AP}
There exist a compact abelian group $G$, an element $\tau$ of $G$, and 
a continuous function $f$ on $G$ such that $a_n=f(\tau^n)$ for every 
$n\in\Z$.
\item \label{item:exponential}
The sequence $\ua$ is a uniform limit of linear combinations 
of exponential sequences.
\item
\label{item:third}
The family of translated sequences $\{\sigma^k\ua\colon k\in\Z\}$ is 
relatively compact under the $\ell^\infty$-norm.
\end{enumerate}
A sequence satisfying these properties is called \emph{almost periodic}.
\end{proposition}

The family of almost periodic sequences forms a sub-algebra of 
$\ell^\infty$, 
closed under complex conjugation, the shift and under uniform limits.
We denote this sub-algebra by $\AP$.

\subsection{$2$-step nilsequences}

\begin{definition}
\label{def:nilsequ}
When $G$ is a group, $G_2$ denotes its commutator subgroup, that is 
the subgroup spanned by elements of the form $ghg\inv h\inv$ for $g,h\in G$.
$G$ is {\em $2$-step nilpotent} if $G_2$ is included in the center of $G$. 

If $G$ is a $2$-step nilpotent Lie group and $\Gamma\subset G$ is 
a discrete and cocompact subgroup, then the compact 
manifold $X = G/\Gamma$ is called a {\em $2$-step nilmanifold}.  
The action of $G$ on $X$ by 
left translation is written $(g,x)\mapsto g\cdot x$ for $g\in G$ and 
$x\in X$.
A $2$-step nilmanifold 
$X=G/\Gamma$, endowed with the translation $T\colon x\mapsto\tau\cdot 
x$ by some fixed element  $\tau$ of $G$ is called a \emph{$2$-step nilsystem}.

If $f\colon X\to\C$ is a continuous function, $\tau\in G$ and 
$x_0\in X$, then $(f(\tau^n\cdot x_0)\colon n\in\Z)$ is a 
{\em basic $2$-step nilsequence}.
A {\em $2$-step nilsequence} 
is a uniform limit of basic $2$-step nilsequences.
\end{definition}

With small modifications, this generalizes condition~\eqref{item:AP} in 
Proposition~\ref{prop:AP}.  
There is a technical difficulty that explains the need to 
take uniform limits of basic nilsequences, rather than just nilsequences.  
Namely, $\AP$ is closed under the uniform norm, 
while the family of basic $2$-step nilsequences is not.
An inverse limit of rotations on compact abelian Lie groups is 
also a rotation on a compact abelian group, but the same result does not 
hold for nilsystems: the inverse limit of a sequence of $2$-step nilmanifolds 
is not, in general, the homogeneous space of some locally compact 
abelian group (see~\cite{R}).
  The other conditions in Proposition~\ref{prop:AP} 
do not have straightforward modifications that generalize for $2$-step 
nilsequences.  After introducing some particular classes of 
$2$-step nilsequences, we develop analogs of these other conditions.

There is an alternate characterization of $2$-step nilsequences:
\begin{proposition*}[see~\cite{HM}]
A sequence $\ua=(a_n\colon n\in\Z)$ is a $2$-step nilsequence if 
there exist an inverse limit $(X,T)$ of $2$-step nilsystems, a 
continuous function $f$ on $X$, and a point $x_0\in X$ such that 
$a_n=f(T^nx_0)$ for every $n\in\Z$.
\end{proposition*}

Definition~\ref{def:nilsequ}, with the obvious changes, defines a $k$-step 
nilsequence.  However, in this article we restrict ourselves 
to studying $2$-step nilsequences. 

\begin{convention} 
We often omit the $2$-step 
from our terminology and just refer to a 
$2$-step nilsequence as a nilsequence.  Similarly, a nilsystem refers 
to a $2$-step nilsystem.
\end{convention}

We introduce some particular $2$-step nilsequences that play a 
particular role in the sequel.

\subsection{Affine systems and quadratic sequences}
\label{subsec:affine}

\begin{notation}
For $t\in\T$, let $\uq(t)=(q_n(t)\colon n\in\Z)$ denote the sequence given by
$$
 q_n(t)=e\bigl(\frac{n(n-1)}2t\bigr)\quad \text{ for } n\in\Z \ .
$$
 This sequence is called a \emph{quadratic exponential sequence}.
\end{notation}
We show that any sequence $\uq(t)$ is a nilsequence.
We first remark that this sequence is almost periodic if and only 
if $t$ is rational and that in this case it is periodic.

Let $G = \Z\times\T\times\T$, $\Gamma = \Z\times \{0\}\times\{0\}$ 
with multiplication defined by 
$$
(m, x, y)\cdot (m', x', y') = (m+m', x+x', y+y'+mx') \ .
$$
Then $G$ is a $2$-step nilpotent Lie group, with 
$G_2=\{0\}\times\{0\}\times\T$, 
 and $\Gamma$ is a discrete cocompact subgroup. 
Set $X=G/\Gamma$.
Fix some $\alpha,\beta\in\T$, and set
$\tau=(1,\alpha,\beta)\in G$ and endow 
$X$ with the translation $T$ by $\tau$.

We can give another description of this system. The map $(m,x,y)\mapsto
(x,y)$ from $G\to \T^2$ induces a homeomorphism from $X$ to $\T^2$. 
Identifying these spaces, $T$ becomes the familiar skew transformation 
$(x,y)\mapsto (x+\alpha,y+x+\beta)$ of $\T^2$.

We describe some nilsequences arising from this system. 
Set $x_0=(0,0)\in X$. 
Let  $f$ be the continuous function on $X$ given by 
$f(x,y)=e(y)$. 
For every $n\in\Z$,
$$
f(T^nx_0)=e\bigl(\frac{n(n-1)}2 
\alpha+n\beta\bigr)=q_n(\alpha)e_n(\beta)\ .
$$
Thus every  product of a quadratic exponential sequence and
an exponential sequence is a basic 
elementary nilsequence
(In particular, an exponential sequence is a basic 
nilsequence).

\subsection{Heisenberg nilsystems}
\label{subsec:Heisenberg}
\subsubsection{Heisenberg groups}
\label{subsec:defHeisenberg}
Let $d\geq 1$ be an integer. The \emph{Heisenberg group of dimension 
$2d+1$} is $G=\R^{2d+1}$, identified with $\R^d\times\R^d\times\R$ and 
endowed with multiplication given by
$$
 (x,y,z)\cdot(x',y',z')=(x+x',y+y',z+z'+\langle x\mid y'\rangle
\quad \text{ where }x,y,x',y'\in\R^d,\ z,z'\in\R
$$
and $\langle\cdot\mid\cdot\rangle$ is the usual inner product on 
$\R^d$. The commutator subgroup $G_2$ of $G$ is 
$\{0\}^d\times\{0\}^d\times\R$ and $G$ is a $2$-step nilpotent Lie 
group.
We set $\Gamma=\Z^{2d+1}$. Then $\Gamma$ is a discrete cocompact 
subgroup of $G$. The nilmanifold $N_d=G/\Gamma$ is called 
the \emph{Heisenberg nilmanifold of dimension $2d+1$}.

We give another representation of this nilmanifold, which is more 
useful in the sequel. 
The subgroup $K=\{0\}^d\times\{0\}^d\times\Z$ 
is normal in  $G$, and is included in $\Gamma$. We thus have that
$N_d=H_d/\Lambda_d$, where $N_d=G/K$, $\Lambda_d=\Gamma/K$. That is,
$$
 H_d=\R^d\times\R^d\times\S^1,\ \Lambda_d=\Z^d\times\Z^d\times\{1\}\ ,
$$
where the multiplication is given by
\begin{equation}
\label{eq:Heisenberg}
(x,y,z)\cdot(x',y',z')=\bigl(x+x',y+y',zz'e(
\langle x\mid y'\rangle)\bigr)
\quad (x,y,x',y'\in\R^d,\ z,z'\in\S^1)\ .
\end{equation}

(Note that  $\S^1$ denotes the circle while 
$\T=\R/\Z$ is called the torus.)
 When $N_d$ is 
endowed with the translation by some element $\tau$ of $H_d$,
$(N_d,T)$ is called a \emph{Heisenberg nilsystem}.   The  
commutator subgroup of $H_d$ is $\{0\}^d\times\{0\}^d\times\S^1$.

\subsubsection{Nilsequences arising from Heisenberg nilsystems}
\label{subsec:nilseqHeisen}

\begin{notation}
For $s,t\in\T$, we set
\begin{equation}
\label{eq:omega}
\kappa(s,t):=\sum_{k\in\Z}\exp(-\pi(t+k)^2) e(ks)\ .
\end{equation}
We write 
$\uomega(\alpha,\beta)=(\omega_n(\alpha, \beta)\colon n\in\Z)$ for the sequence defined by
$$
\omega_n(\alpha,\beta)=
\kappa(n\alpha,n\beta)e(\frac{n(n-1)}2\alpha\beta)\quad \text{ for }n\in\Z \ .
$$
\end{notation}

We show that products of sequences of this type are
are $2$-step nilsequences. 
\begin{remark*}
The function $\kappa$ is linked to the classical theta function
$$
\theta(u,z) = \sum_{n\in\Z}\exp(\pi i z n^2+2\pi inu) 
$$
by
\begin{equation}
\label{eq:omegatheta}
\kappa(s,t) = \exp(-\pi t^2)\theta(s+it, i) \ .
\end{equation}
This is a particular case of well known relations between 
Heisenberg groups and theta functions.
\end{remark*}

We start by explaining how to 
build explicit continuous functions on $N_d$.  It is 
not as easy to do this as it is for affine systems, 
because the nilmanifold $H_d$ is 
not homeomorphic to a torus.

Let $\phi$ be a continuous function on $\R^d$, tending to zero 
sufficiently fast at infinity. Define a function $\tilde f$ on 
$H_d$ by
$$
 \tilde f(x,y,z)=z\sum_{k\in\Z^d}\phi(y+k)e(\langle k\mid 
 x\rangle)\quad
(x,y\in\R^d,\ z\in\S^1)\ .
$$
We note that for all $k,\ell\in\Z^d$, we have 
$\tilde f\bigl((x,y,z)\cdot(k,\ell,1)\bigr)=\tilde f(x,y,z)$.
In other words, the function $\tilde f$ on $G$ is invariant under 
right translations by elements of $\Lambda_d$, and thus it induces a 
continuous function $f$ on $N_d$.

The most interesting case is when $\phi$ is the Gaussian function 
$\phi(x)=\exp(-\pi\norm x^2)$, where $\norm\cdot$ is the usual 
euclidean norm in $\R^d$.

In this case we have:
\begin{equation}
\label{eq:fomega}
 \tilde f(x,y,z)=z\prod_{j=1}^d\kappa(x_j,y_j)\ .
\end{equation}

Set $\tau=(\alpha,\beta,e(\gamma))\in H_d$, where 
$\alpha=(\alpha_1,\dots,\alpha_d)$ and 
$\beta=(\beta_1,\dots,\beta_d)$ belong to $\R^d$ and $\gamma\in\T$
and define $T$ to be 
the translation by $\tau$ on $N_d$. For every $n\in\Z$, we have 
$$
\tau^n=(n\alpha,n\beta, e(n\gamma)e\bigl(\frac 
{n(n-1)}2\langle\alpha\mid\beta\rangle\bigr))\ .
$$
Let $x_0$ be the image 
in $N_d$ of the element $(0,0,1)$ of $H_d$ and $f$ the function on 
$N_d$ associated to the function $\tilde f$ defined 
by~\eqref{eq:fomega}.
Then for every $n\in\Z$, we have
\begin{multline}
\label{eq:sequHeisenberg}
 f(T^n x_0)=e(n\gamma)e\bigl(\frac 
{n(n-1)}2\langle\alpha\mid\beta\rangle\bigr)
\prod_{j=1}^d\kappa(n\alpha_j,n\beta_j)\\
=e_n(\gamma)
\omega_n(\alpha_1,\beta_1).\cdots.\omega_n(\alpha_d,\beta_d)\ .
\end{multline}

\subsection{Some counterexamples}
There are other natural quadratic-type sequences 
that are particularly simple to construct, and so it 
is natural to ask if these are (elementary) nilsequences.
For example, consider the sequence 
$(e(\lfloor n\alpha\rfloor n\beta)\colon n\in\Z)$.  
This is not a nilsequence, as it can not be
obtained as a uniform limit of sequences associated to continuous functions
on nilmanifolds.  On the other hand, this 
sequence is a ``Besicovitch'' elementary nilsequence, meaning that it is a
limit  of elementary nilsequences in quadratic 
averages.  
This can be seen by considering 
the simpler sequence $(e(\lfloor n\alpha\rfloor  \beta)\colon n\in\Z)$.  
Although at 
first glance it looks like an almost periodic sequence,
it is not: under the uniform norm, 
the sequence  $(e(\lfloor n\alpha +\gamma\rfloor \beta)\colon n\in\Z)$ 
does not depend continuously on $\gamma$.  
The family of shifted sequences $(a_{n+k}\colon n\in\Z)$ is
not relatively compact under the uniform norm, but only 
in quadratic average.
A similar problem occurs with the 
sequence $(e([n\alpha]n\beta)\colon n\in\Z)$.
A more complete explanation is given in Appendix~\ref{ap:non-example}.

Objects similar to our 
nilsequences arise in several other contexts.  For example, 
Green and Tao~\cite{GT} use nilsequences to give 
asymptotics for the number of progressions of length $4$ in the 
primes.  Their definition is slightly different:
the underlying group $G$ is assumed to be connected and 
simply connected.  Moreover, in all their uses 
of nilsequences the function defining it is taken to 
be Riemann integrable.   
This weaker hypothesis on the function (we assume continuity)
accounts for the difference in sequences that arise.  In particular, 
sequences with integer parts can be obtained from Riemann integrable 
functions, but not from continuous functions.  
Perhaps more important is the difference in point of view.  
In Green and Tao, they consider finite (albeit long) sequences taken from 
a nilsequence.  The importance is that this point of view 
is local and in this context, quadratic sequences are the natural local 
model for $2$-step nilsequences.  
In giving a global model for $2$-step nilsequences, this program 
no longer can be carried through. 

In~\cite{BL}, Bergelson and Leibman 
show that all generalized polynomials arise by 
evaluating a piecewise polynomial mapping on a nilmanifold.
However, in order to obtain all generalized polynomials, 
they necessarily must allow functions 
with discontinuities.  Our nilsequences are recurrent, 
while some generalized polynomials take on non-recurring values.

\section{Results}
\label{sec:results}
\subsection{First properties}
\label{subsec:first}
\begin{notation} We 
denote the  family by of ($2$-step) nilsequences by 
$\Nil$. 
\end{notation}

Some classical properties of $2$-step nilsystems are recalled in 
Section~\ref{sec:basics}. 

\begin{proposition}
\label{prop:closed}\strut
\begin{itemize}
\item
The family of nilsequences is a 
sub-algebra of $\ell^\infty(\Z)$ that contains $\AP$ and is invariant under 
complex conjugation, shift and uniform limits. 
\item
If $\ua$ is a nilsequence and $(I_j)_{j\in\N}$ is 
a sequence of intervals of $\Z$ whose lengths tend to infinity,
then the limit
$$
 \Av(\ua):=\lim_{j\to\infty}\frac 1{|I_j|}\sum_{n\in I_j}a_n
$$
exists and does not depend on choice 
of the sequence $(I_j)$ of intervals. 
\end{itemize}
\end{proposition}

If $\ua$ and $\ub$ are nilsequences, then the product sequence 
$\ua\,\overline{\ub}:=(a_n\overline{b_n}\colon n\in\Z\}$ is a 
nilsequence.
This, together with the second part of Proposition~\ref{prop:closed}
justifies the following definition:
\begin{definition}
\label{def:quadraticnorm}
The \emph{inner product} of two nilsequences $\ua$ and $\ub$ is
defined to be
$$
 \langle\ua\mid\ub\rangle:=\Av(\ua\,\overline{\ub})=
\lim_{N\to+\infty}\frac{1}{N}\sum_{n=0}^{N-1}a_n\overline{b_n}\ .
$$
The \emph{quadratic norm} of the nilsequence $\ua$ is
$$
 \norm{\ua}_2:=\bigl(\Av(|\ua|^2)\bigr)^{1/2}=
\lim_{N\to+\infty}\frac{1}{N}\Bigl(\sum_{n=0}^{N-1}|a_n|^2\Bigr)^{1/2} \ .
$$
\end{definition}
In Section~\ref{subsec:quadratic-norm} we check that $\norm\cdot_2$ is 
a norm on the space $\Nil$ of nilsequences.

\subsection{Elementary $2$-step nilsequences}
There is no direct generalization of the second 
property of Proposition~\ref{prop:AP} for $2$-step nilsequences, 
as there is no strict
analog of exponential sequences.  However, some particular 
nilsequences play a similar role, and these are 
the building blocks for all nilsequences.

\begin{definition}
\label{def:elementary}
A $2$-step nilmanifold $X$ is an {\em elementary ($2$-step) nilmanifold}
if it can be written as $X=G/\Gamma$ as in 
Definition~\ref{def:nilsequ}, where $G_2$ is the circle group $\S^1$.
An \emph{elementary nilsystem} is a  nilsystem $(X=G/\Gamma,T)$ such 
that $X$ is an elementary  nilmanifold.

When $X=G/\Gamma$ is an elementary nilmanifold (written such 
that $G_2 = \S^1$), we write $\CC_1(X)$ 
for the family of continuous functions $f$ on $X$ such that
\begin{equation}
\label{eq:C1}
f(u\cdot x) = uf(x) \text{ for all } u\in G_2=\S^1\text{ and all } x\in X \ .
\end{equation}

If $(X,T)$ is an elementary nilsystem, 
$f\in\CC_1(X)$  and $x_0\in X$, the 
sequence $\bigl(f(T^n x_0)\colon n\in\Z\bigr)$ is a 
\emph{basic elementary nilsequence}.
A uniform limit of basic elementary nilsequences is an
\emph{elementary nilsequence}.

By convention, every almost periodic sequence is also considered to be
an elementary nilsequence.
\end{definition}

There are several things left imprecise in the definition of an elementary 
nilmanifold. 
First, the same nilmanifold $X$ can be written as 
$G/\Gamma$ in different ways, and the term ``elementary nilmanifold'' 
refers in fact to a given presentation of $X$ and not only to the 
nilmanifold. 
Moreover, when we say that ``$G_2=\S^1$'' we do not mean only that 
these groups are isomorphic, but also that we make the choice of a 
particular isomorphism that defines the identification.

\begin{definition}
Let $(X=G/\Gamma,T)$ be an elementary nilsystem. The same system but 
with the opposite 
identification of $G_2$ with $\S^1$ is called the \emph{conjugate system} of $X$ and is 
written $(\bar X,T)$.
We note that for $f\in\CC_1(X)$, we have that $\bar f\in\CC_1(\bar 
X)$.
\end{definition}

Under the usual identification of $\T$ with $\S^1$, the affine 
system $(X,T)$ introduced in Section~\ref{subsec:affine} is an elementary 
nilsystem, the function  $f$ belongs to $\CC_1(X)$ and the sequence
$(q_n(\alpha)e_n(\beta)\colon n\in\Z)$ is a basic elementary 
nilsequence.

The Heisenberg system $(N_d,T)$ introduced in 
Section~\ref{subsec:defHeisenberg} is an elementary nilsystem, and 
the sequence defined in~\eqref{eq:sequHeisenberg} is a basic 
elementary nilsequence.

\begin{proposition}
\label{prop:elementary-closed}
The family of elementary nilsequences is closed under complex 
conjugation, the shift, taking products and uniform limits.
\end{proposition}
All these properties follow immediately from the 
definitions, other than the closure under products.  This is proved in 
Section~\ref{sec:products} (Proposition~\ref{prop:equalortho2}).

The next Proposition (proved in Section~\ref{subsec:decomposition}) 
is a nilsequence version of the characterization 
of almost periodic sequences given by Part~\eqref{item:exponential} of 
Proposition~\ref{prop:AP}:

\begin{proposition}
\label{prop:decomposition}
Every nilsequence is a uniform limit of finite sums of 
elementary nilsequences.
\end{proposition}

In this statement, we can obviously substitute ``basic 
elementary nilsequences'' for elementary nilsequences.

The family of $2$-step nilsequences has no 
characterization similar to
property~\eqref{item:third} of Proposition~\ref{prop:AP}.  However, 
the smaller class of elementary nilsequences does (see 
Section~\ref{sec:charactelemen2}
for the proof):

\begin{theorem}
\label{th:charactelemen}
Let $\ua=(a_n\colon n\in\Z)$ be a bounded sequence.  
The following are equivalent:
\begin{enumerate}
\item
\label{it:elementary} $\ua$ is an elementary nilsequence.

\item
\label{it:compact} There exists a compact 
(in the norm topology) subset $K\subset 
\ell^\infty(\Z)$ such that for all $k\in \Z$, there exists $t\in T$ 
such that the sequence $\ue(t)\sigma^k\ua=\bigl(e(nt)a_{n+k}\colon 
n\in\Z\bigr)$  belongs to $K$.
\end{enumerate}
\end{theorem}

\subsection{The class of an elementary nilsequence}

Theorem~\ref{th:charactelemen} makes it clear that 
the fundamental objects are not elementary nilsequences, 
but rather elementary nilsequences considered relative to 
almost periodic sequences.  This 
motivates the next definition.
\begin{definition}
Let $\ua= (a_n\colon n\in\Z)$ be a bounded sequence. The \emph{class} 
$\CS(\ua)$ of $\ua$ is the norm closure in $\ell^\infty$ of the vector space
spanned by sequences of the form
 $\ue(t)\sigma^k\ua=\bigl(e(nt)a_{n+k}\colon n\in\Z\bigr)$ for 
 $t\in\T$ and $k\in\Z$.
\end{definition}

The class of an elementary nilsequence contains only elementary 
nilsequences. If $\ua$ is a non-identically zero almost periodic 
sequence, then $\CS(\ua)=\AP$.
\begin{remark*} The reason for the introduction of the shift in the 
definition of the class $\CS(\ua)$ is not obvious at this point but 
is clarified later (proof of Corollary~\ref{cor:NXirreducible}).
 In fact, if an elementary sequence $\ua$ is 
such that the sequence $(|a_n|\colon n\in\Z)$ is bounded from below, then the shift 
can be removed from the definition:
$\CS(\ua)$ is spanned by sequences of the form $\ue(t)\ua$.  
However, there are very few elementary nilsequences of this 
type: they all are products of an almost periodic sequence and a 
quadratic sequence $\uq(t)$, as defined in Section~\ref{subsec:affine}.
The situation is different if we use the quadratic norm (see the proof 
of Theorem~\ref{th:span2} in Section~\ref{subsec:proofspan2}).
\end{remark*}

In Section~\ref{subsec:proofequalortho} 
we show:

\begin{theorem}
\label{th:equalortho}\strut
\begin{enumerate}
\item\label{it:equalortho}
The classes of two elementary nilsequences are either identical or are
orthogonal. This means that if $\ua$ and $\ub$ are non-identically zero
elementary 
nilsequences, then:
\begin{itemize}
\item
If $\ub\in\CS(\ua)$, then $\CS(\ub)=\CS(\ua)$.

\item If $\ub\notin\CS(\ua)$, then $\Av(\ua\overline{\ub})=0$, that is
$$
 \frac 1N\sum_{n=0}^{N-1}a_n\overline{b_n}\to 0\text{ as }N\to\infty\ .
$$
\end{itemize}
\item\label{it:open}
The classes form a partition of the set of elementary nilsequences into 
open subsets.
\item\label{it:products}
If $\ua$ and $\ub$ are elementary nilsequences, then there exists an 
elementary nilsequence $\uc$ such that 
$\CS(\ua)\cdot\CS(\ub):=\{\ua'\ub'\colon\ua'\in\CS(\ua),\ 
\ub'\in\CS(\ub)\}$ is included in $\CS(\uc)$.
\item\label{it:groups}
Endowed with this multiplication, the family of classes is an abelian 
group.
\end{enumerate}
\end{theorem}

\subsection{Classification}
\label{subsec:classific}
In light of the 
 last two parts of Theorem~\ref{th:equalortho}, it would be nice 
to classify elementary nilsequences.
Optimally, one would like to classify the elementary nilsequences 
by giving an explicit 
representative in every class.  
Unfortunately, this is not possible as there is no 
reasonable (i.e. Borel) way to choose this representative. 
However, we can partially succeed in giving a classification, 
by reducing to some particular nilsequences.

Our classification of elementary nilsequences is related 
to, but not identical to, a classification of elementary nilsystems. 
In this classification, non-isomorphic 
nilsystems may give rise to the same class of 
nilsequences.

\begin{notation}
 Let $d\geq 1$ be an integer. We write $J_{2d}$ for the 
$2d\times 2d$ matrix 
$$
 J_{2d}=\begin{pmatrix}
0 & \id_d\\ -\id_d & 0
\end{pmatrix}
$$
where the $0$ represents the $d\times d$ zero matrix and $\id_d$ the 
$d\times d$ identity matrix.
 
Let $\Sp_{2d}(\Q)$ denote the group of $2d\times 2d$ matrices $M$ with 
rational entries which preserve the antisymmetric bilinear form 
associated to $J_{2d}$, that is, with
$$
 M^tJ_{2d}M=J_{2d}\ .
$$
In an equivalent way, $M$ can be written in $d\times d$ blocks as
$$
 M=\begin{pmatrix}
A & B \\ C & D
\end{pmatrix}
$$
where  $A^tC$ and $D^tB$  are symmetric and 
$D^tA-B^tC = \id_d$.
\end{notation}
\begin{convention}
In the next theorem and in the sequel, the empty product is by 
convention set to be $1$.
\end{convention}

\begin{theorem}
\label{th:classification}

Let $t\in\T$ be either equal to $0$ or irrational, let
$d\geq 0$ be an integer  and let  $\alpha_1,\dots,\alpha_d, \beta_1, \ldots, 
\beta_d$ be reals, rationally  independent  modulo $1$.
Then the sequence 
$\uq(t)\,\uomega(\alpha_1,\beta_1).\cdots.\uomega(\alpha_d,\beta_d)$
is an elementary nilsequence.

Every elementary nilsequence belongs to the class of a 
sequence of this type.

Let $t,t'\in\T$, $d,d'\geq 0$ 
and $\alpha_1,\dots,\alpha_d,\beta_1,\dots,\beta_d$ and
$\alpha'_1,\dots,\alpha'_{d'},\beta'_1,\dots,\beta'_{d'}$ be 
rationally independent modulo $1$. Then the sequences 
$\uq(t)\,\uomega(\alpha_1,\beta_1).\cdots.\uomega(\alpha_d,\beta_d)$ and 
$\uq(t')\,\uomega(\alpha'_1,\beta'_1).\cdots.\uomega(\alpha'_{d'},\beta'_{d'})$ 
belong to the same class if and only if $d=d'$ and 
there exist
a $2d\times 2d$ matrix $Q\in\Sp_{2d}(\Q)$
and integers $m,k_1,\dots,k_d,\ell_1,\dots,\ell_d$ with $m\geq 1$
such that 
$$
 Q\begin{pmatrix}
\alpha_1+k_1/m\\ \dots \\ \alpha_d+k_d/m\\ \beta_1+\ell_1/m
\\ \beta_d+\ell_d/m
\end{pmatrix}=\begin{pmatrix}
\alpha'_1\\ \dots \\\alpha'_d\\ \beta'_1\\ \dots \\ \beta'_d
\end{pmatrix}
\text{ and }
m(t-t')= \sum_{i=1}^d(k_i\beta_i-\ell_i\alpha_i)\bmod 1\ .
$$
\end{theorem}

If $d=0$ in the first part of the theorem, then the sequence is $\uq(t)$ 
and the hypothesis of 
independence is vacuous. If $d=d'=0$ in the second part of the theorem, 
then  the sequences are $\uq(t)$ and $\uq(t')$ and the condition means simply that $t-t'$ is rational.

\subsection{Density results}

Together, the examples of nilsequences introduced in 
Sections~\ref{subsec:affine} ad~\ref{subsec:nilseqHeisen} span the $2$-step nilsequences 
and so we give them a name:

\begin{notation}
We write $\CM$ for the family of sequences of the form
$$
\ue(s)\uq(t)\uomega(\alpha_1,\beta_1).\cdots .
\uomega(\alpha_d,\beta_d)
$$
where $s,t\in\T$, $d\geq 0$ is an integer and 
$\alpha_1,\dots,\alpha_d,\beta_1,\dots,\beta_d$ are reals that are rationally 
independent modulo $1$.
\end{notation}

\begin{theorem}
\label{th:span}
The space of $2$-step nilsequences is the closed shift invariant 
linear space spanned by the family $\CM$ of sequences.
\end{theorem}

For the quadratic norm, 
we have a density result analog of Theorem~\ref{th:span}, without 
needing to include shifts of the sequences:
\begin{theorem}
\label{th:span2}
The linear span of the family $\CM$ of sequences in dense in the space 
$\Nil$ of $2$-step nilsequences under the quadratic norm.
\end{theorem}

However, Theorem~\ref{th:classification} shows that $\CM$ is not a basis for $\Nil$.

\subsection{Higher order nilsequences}

Some of the analysis in this paper for $2$-step nilsequences 
can be carried out for higher levels.  
We can define elementary $k$-step nilsequences by simply 
replacing $G_2$ by $G_k$.  Then arbitrary 
$k$-step nilsequences can be approximated by elementary ones. 
However, the problem then becomes finding the analog of the 
span of (elementary) $k$-step nilsequences and 
giving a characterization of (elementary) $k$-step nilsequences. 

One of the advantages of working with $2$-step nilsequences is 
that we have the family of theta functions 
and associated theta sequences that simplify the classification. 
Representations of $k$-step nilpotent groups for $k\geq 3$ are 
significantly more complicated, and so the related 
classification of $k$-step nilsequences is a difficult problem.

\section{Basic properties of nilsystems and nilsequences}
\label{sec:basics}
\subsection{Nilmanifolds and nilsystems}

We recall some classical properties of nilmanifolds and nilsystems. 
For nilmanifolds, we refer to the seminal paper by Malcev~\cite{M}. 
For nilsystems, we refer to~\cite{AGH}, \cite{Pa1}, \cite{Pa2} and 
\cite{Les}; these 
papers deal only with the case of connected nilpotent Lie groups and 
their results were extended to the non-connected case in \cite{Lei2}.

In section~\ref{subsec:affine}, we saw that every exponential 
sequence is a nilsequence. Therefore, $\Nil\supset\AP$.
The family of $2$-step nilsystems is 
closed under Cartesian products.  This, together with 
Definition~\ref{def:nilsequ}, immediately implies the first statement 
of Proposition~\ref{prop:closed}.

In the rest of this section, 
 $(X=G/\Gamma,T)$ is assumed to be a $2$-step nilsystem, where $T$ is the left 
translation by $\tau\in G$ and $x_0$ is a point of $X$.

Let $Y$ be the closed orbit of $x_0$, that is the closure in 
$X$ of $\{T^nx_0\colon n\in \Z\}$. Then $(Y,T)$ can be given the 
structure of a $2$-step nilsystem or of a rotation on a compact 
abelian group.   
More precisely, there exists a 
closed subgroup $H$ of $G$ containing $\tau$ such that 
$\Lambda:=\Gamma\cap H$ is cocompact in $H$ and $Y$ can be identified 
with $H/\Lambda$. $H$ is a Lie group, either abelian or $2$-step 
nilpotent. 

By substituting $H$ for $G$ and $\Lambda$ for 
$\Gamma$ in the definitions of basic $2$-step nilsequences and of 
elementary nilsequences, we can restrict to the case that the orbit 
$\{T^nx_0\colon n\in\Z\}$ of $x_0$ is dense in $X$. It is classical 
 that this hypothesis implies:
\begin{enumerate}
\item\label{it:minimal}
$(X,T)$ is minimal (that is, 
every orbit is dense) and uniquely ergodic: the unique $T$-invariant 
probability measure is the \emph{Haar measure} $\mu$ of $X$, which is 
invariant under the action of $G$ on $X$.
\end{enumerate}
Henceforth we assume that our nilsystem 
satisfies this assumption.

This immediately implies 
the second statement of Proposition~\ref{prop:closed}:  if a basic 
nilsequence $\ua$ is defined as in Definition~\ref{def:nilsequ},
then 
$$
 \Av(\ua)=\int f\,d\mu
$$
and the existence of the average in the case of a general nilsequence 
follows by density.
\subsection{The quadratic norm}
\label{subsec:quadratic-norm}
The quadratic norm $\norm\ua_2$ of a $2$-step nilsequence was 
introduced in Section~\ref{subsec:first} 
(Definition~\ref{def:quadraticnorm}); we check here 
that it is a norm.

We only need to show that if $\ua$ is a non-identically zero 
nilsequence, then $\norm{\ua}_2\neq 0$.
Let $\ua$ be a non-identically zero basic nilsequence. 
Write it as $\bigl(f(T^nx_0)\colon 
n\in\Z\bigr)$, where $(X,T)$ is a minimal nilsystem, $f\in\CC(X)$ and 
$x_0\in X$.
By minimality, we have that $\norm{\ua}_\infty=\norm f_\infty$ and 
that 
$$
\text{for every }\delta \text{ with }0<\delta<\norm{\ua}_\infty
\text{,  the set }
\{n\in\Z\colon |a_n|\geq\delta\}\text{ is syndetic.}
$$
This property is stable under uniform limits and thus holds for all
nilsequences.  The announced result follows.\qed

We remark that $\norm\ua_2\leq\norm\ua_\infty$ for every 
nilsequence $\ua$.

\subsection{Reduced form of a $2$-step nilsystem}
We can make some further simplifying assumptions 
about the nilsystem, as in~\cite{BHK}. 
We only sketch the justification of 
these reductions.

Let $G_0$ be the connected component of the identity in $G$ and let $G_1$ 
be the subgroup of $G$ spanned by $G_0$ and $\tau$. Since $G_1$ is open and 
the projection $G\to G/\Gamma=X$ is an open map, the image of $G_1$ 
under this projection is an open subset of $X$. This subset is 
invariant under $T$ and thus is equal to $X$ by minimality. Therefore, 
substituting $G_1$ for $G$ and $\Gamma\cap G_1$ for $\Gamma$, we can 
reduce to the case that
\begin{enumerate}
\setcounter{enumi}{1}
\item
\label{item:identity}
$G$ is spanned by the connected component $G_0$ of the identity 
 and $\tau$.
\end{enumerate}

Let $\CZ(G)$ be the center of $G$. 
We can substitute $G/(\Gamma\cap\CZ(G))$ for $G$,
and $\Gamma/(\Gamma\cap\CZ(G))$ for $\Gamma$. 
Therefore we can reduce to the case that
\begin{enumerate}
\setcounter{enumi}{2}
\item
\label{item:trivial-center}The intersection of $\Gamma$ with the 
center of $G$ is trivial.
\end{enumerate}
This implies that $\Gamma$ does not contain any nontrivial normal
subgroup of $G$ and thus that $G$ can be viewed as a group of transformations of
$X$.

These properties imply that
\begin{enumerate}
\setcounter{enumi}{3}
\item
$\Gamma$ is abelian.
\item
$G_2$ is a finite dimensional torus.
\end{enumerate}

\begin{definition}
Let $(X=G/\Gamma,T)$ be a ($2$-step) nilsystem $(X,T)$, 
where $T$ is the 
translation by $\tau\in G$. 
We say that $(X,T)$ is \emph{in reduced form} if 
$G,\Gamma,\tau$ satisfy the hypotheses~\eqref{item:identity} 
and~\eqref{item:trivial-center}.
\end{definition}

By the preceding reductions, we see that in the definition of $2$-step 
nilsequences we can restrict to the case that the nilsystem is minimal 
and in reduced form.
We claim that 
 in the definition of an elementary nilsequences we 
can restrict ourselves to minimal elementary nilsystems written in 
reduced form.

 Indeed, let 
$(X=G/\Gamma,T)$ be an elementary nilsystem, $x_0\in X$,
$f\in\CC_1(X)$ and let $\ua=(f(T^nx)\colon n\in\Z)$.
We can assume without loss that this sequence is not almost periodic.

Let $Y$ be the closed orbit of $x_0$ under $T$. 
Then $(Y,T)$ is a not a rotation on a compact abelian group and thus 
is a $2$-step nilsystem. We write $Y=G'/\Gamma'$, where $G'$ is a 
closed subgroup of $G$ and $\Gamma'=\Gamma\cap G'$.
As this system is minimal, we can make the 
second reduction and assume that property~\eqref{item:identity} 
holds. 
Then  $G'_2$ is a nontrivial closed connected subgroup of $G_2=\S^1$, 
and thus $G'_2=G_2$ and $(Y,T)$ is an elementary nilsystem.
Moreover, the restriction of $f$ to $Y$ belongs to $\CC_1(Y)$. 

Now we make the last reduction, and write $Y=G''/\Gamma''$ by taking 
the quotient of $G'$ and $\Gamma'$ by $\CZ(G')\cap\Gamma'$.
The existence of the non-identically zero function 
$f$ in $\CC_1(Y)$ implies that $G'_2\cap \Gamma'$ is trivial, and thus
 the natural projection $G'_2\to G''_2$ is an isomorphism.
We thus have that $G''_2$ is the circle group and that $f$ belongs 
to $\CC_1(Y=G''/\Gamma'')$. Our claim is proved.

\begin{convention}
Henceforth, all nilsystems are implicitly assumed to be minimal and 
in reduced form.
\end{convention}

The {\em Kronecker factor} of an ergodic system 
$(X,T,\mu)$ is the factor whose $\sigma$-algebra is generated 
by all the eigenfunctions of $T$.  
We note the following result for later use (see for 
example~\cite{Lei2}).
\begin{proposition}
\label{prop:cond-min}
Let
$X=G/\Gamma$ where $G$ is a 
$2$-step nilpotent Lie group, $\Gamma$ is a closed cocompact subgroup 
and let $T$ be translation by some $\tau\in G$. Let $X$ be endowed 
with its Haar measure $\mu$. Assume that 
property~\eqref{item:identity} is satisfied.  Then the following 
properties are equivalent:
\begin{enumerate}
\renewcommand{\theenumi}{\arabic{enumi}}
\item
$(X,T)$ is minimal.
\item
The subgroup of $G$ spanned by $\tau$ and $\Gamma$ is dense in $G$.
\item
The subgroup of $G$ spanned by $\tau$ and $\Gamma G_2$ is dense in $G$.
\end{enumerate}
If these properties are satisfied, then 
the Kronecker factor of the ergodic system $(X,\mu,T)$ is the compact 
abelian group $Z=G/G_2\Gamma$ endowed with translation by the 
image of $\tau$.
\end{proposition}

\subsection{The examples of Section~\ref{subsec:affine} 
and~\ref{subsec:Heisenberg}, revisited}
\label{subsec:backexamples}
By the preceding criteria, the affine system defined in Section~ 
\ref{subsec:affine} is minimal if and only if $\alpha$ is irrational. 
We remark that this system is written in reduced form.

The first presentation of the Heisenberg system $N_d=G/\Gamma$ given 
in Section~\ref{subsec:Heisenberg} was not reduced. The reduction to 
the presentation   $N_d=H_d/\Lambda_d$ we carried out in that Section is 
exactly the reduction explained above to the case that 
hypothesis~\eqref{item:trivial-center} holds. In this presentation, 
$N_d$ is written in reduced form.
By the criterion of Proposition~\ref{prop:cond-min}, 
the nilsystem $(N_d,T)$ is minimal if and only 
if the reals $\alpha_1,\dots,\alpha_d,\beta_1,\dots,\beta_d$ are rationally 
independent modulo $1$.

\subsection{Decomposition into elementary nilsequences}
\label{subsec:decomposition}
In this Section, we prove Proposition~\ref{prop:decomposition} and 
introduce some notation used in the sequel.

Let $(X=G/\Gamma,T)$ be a minimal $2$-step nilsystem in reduced 
form.

The restriction of the left translation $G\times X\to X$ to $G_2$ 
defines a continuous action of $G_2$ on $X$, and this 
action commutes with $T$.
The quotient of $X$ under this action is the compact abelian group 
$Z=G/G_2\Gamma$.

For every character $\chi$ of $G_2$, let $\CC_\chi(X)$ denote the space 
of continuous functions $f$ on $X$ satisfying
\begin{equation}
\label{eq:Cchi}
f(u\cdot x)=\chi(u)\,f(x)\quad(u\in G_2,\ x\in X)\ .
\end{equation}
If $X$ is an elementary nilsystem, that is, if $G_2=\S^1$, we 
identify the dual group $\wh{G_2}$ with $\Z$ and the 
notation in~\eqref{eq:Cchi} coincides with the notation in~\eqref{eq:C1}
of Definition~\ref{def:elementary}.

Note that  $\CC_0(X)$ is the space of functions on $X$ which 
factorize through $Z$.   Every space $\CC_\chi(X)$ is invariant under 
$T$ and under multiplication by functions belonging to $\CC_0(X)$. 
The spaces $\CC_\chi(X)$ corresponding to different characters 
$\chi$ are orthogonal with respect to the inner product of 
$L^2(\mu)$.

Let $f\in\CC(X)$.  For every 
$\chi\in\wh{G_2}$, we define a function $f_\chi$ on $X$ by
$$
 f_\chi(x)=\int_{G_2}f(u\cdot x)\,\overline{\chi(u)}\,dm_{G_2}(u)\ ,
$$
where $m_{G_2}$ denotes the Haar measure of $G_2$. Clearly, for every 
$\chi\in\wh{G_2}$. the function $f_\chi$ belongs to the space 
$\CC_\chi(X)$. Moreover,
$$
 f=\sum_{\chi\in\wh{G_2}}f_\chi
$$
where the series converges in $L^2(X)$, and in uniform norm 
whenever $f$ is sufficiently smooth. 
By density, every continuous function 
$f$ on $X$ can be approximated uniformly by finite sums of the 
form
$$
\sum_{i=1}^k f_i\ , \text{ where }
f_i\in\CC_{\chi_i}(X) \text{ and }\chi_i\in\wh{G_2}
\text{ for every }i\ .
$$
Therefore, every basic nilsequence $\ua$ arising from $X$ as 
in Definition~\ref{def:nilsequ} can be approximated 
uniformly by finite sums of nilsequences of the form
$\ub= (f(T^nx_0)\colon n\in\Z)$, where $f$ belongs to $\CC_\chi(X)$ for some 
$\chi\in\wh{G_2}$. 

We are left with describing 
sequences $\ub$ of this form. If $\chi$ is trivial,
then this sequence is almost periodic. Assume that $\chi$ is not 
trivial. 
Let  $Y$ be the 
quotient of $X$ under the action of the subgroup $\ker(\chi)$ of 
$G_2$, $p\colon X\to Y$ be the natural projection and $S$ be the 
transformation induced by $T$ on $Y$. Then it is immediate that 
$(Y,S)$ is an elementary nilsystem. Moreover, every 
$f\in\CC_\chi(X)$ can be written as $f=\phi\circ p$ for some 
continuous function $\phi$ on $Y$ belonging to $\CC_1(Y)$.
Therefore, the sequence $\ub$ is equal to
$(\phi(S^ny_0))$, where $y_0=p(x_0)$.
This completes 
the proof of Proposition~\ref{prop:decomposition}.\qed

\section{The Bohr extension of an elementary nilsystem}
\label{sec:universal}

\subsection{The space of nilsequences associated to an elementary nilsystem}
\strut

\begin{definition} 
Let $(X=G/\Gamma,T)$ be an elementary nilsystem (minimal and written 
in reduced form).  We write $\CN^2(X)$ for the closed linear subspace 
of $\ell^\infty(\Z)$ spanned by sequences of the form 
$\bigl(f(T^nx_0)\colon n\in\Z\bigr)$ where 
$f\in\CC_1(X)$ and $x_0\in X$.
\end{definition}

\begin{lemma}
\label{lem:basepoint}
Let $(X=G/\Gamma,T)$ be an elementary nilsystem. 
\begin{enumerate}
\item\label{it:invarNilX}
 $\Nil(X)$ is 
invariant under the shift and multiplication by almost periodic sequences.
\item\label{it:basepoint}
For every $x_1\in X$, 
$\CN^2(X)$ is the closed linear subspace 
of $\ell^\infty(\Z)$ spanned by sequences of the form 
$\bigl(e(nt)f(T^nx_1)\colon n\in\Z\bigr)$ where $t\in\T$ and
$f\in\CC_1(X)$.
\end{enumerate}
\end{lemma}
\begin{proof}\eqref{it:invarNilX}
The invariance under the shift is obvious. 
It remains to show that every sequence  $\ua$  of the form 
$(e(nt)f(T^nx_0)\colon n\in\Z)$ 
for some $f\in\CC_1(X)$ and $t\in\T$ belongs to $\Nil(X)$.

Set  $K = \{[\tau, g]\colon g\in G\}$ and first we show that $K=\S^1$.
We have that 
$K$ is a subgroup of $\T$.  Since $G$ is spanned 
by $G_0$ and $\tau$, we have that $K$ is connected and 
so is either equal to $\T$ or is trivial.  If $K$ is trivial, 
then $\tau$ belongs to the center $\CZ(G)$ of $G$.  
Given $x\in X$, its closed orbit is thus included in $\CZ(G)\cdot x$, 
and, since it is equal to $X$ by minimality,  we have 
$\CZ(G)\cdot x=X$.
This means that the system $(X,T)$ is a rotation on a compact 
abelian group, a contradiction.  Thus $K=\S^1$.  

Now let $f, x_0,t$ and $\ua$ be as above.
 Pick $g\in G$ such that 
$e(t)=[g,\tau]$ and set $x_1=g\cdot x_0$ and $h(x)=f(g\inv\cdot x)$ for 
every $x\in X$. For every $n\in\Z$,
$$
 h(T^nx_1)=h(\tau^ng\cdot x_0)=h([\tau,g]^ng\tau^n\cdot x_0)=
e(nt)h(g\tau^n\cdot x_0)=e(nt)f(\tau^n\cdot x_0)=a_n
$$
and we are done.

\eqref{it:basepoint}
Fix $x_1\in X$.
Let $M$ denote the closed span of sequences 
of the form $\bigl(e(nt)f(T^nx_1)\colon n\in\Z\bigr)$
 where $t\in\T$ and $f\in\CC_1(X)$.   
We need to show that for every  $x_0\in X$ and every $h\in\CC_1(X)$, we 
have that the sequence  $\bigl(h(T^nx_0)\colon n\in\Z\bigr)$ belongs 
to $M$.

Choose $g\in G$ such that $g\cdot x_1 = x_0$.  
Let $t\in\T$ be such that $e(t)=[\tau, g]\in G_2$ and define 
$f(x) = h(g\cdot x)$. 
Then for $n\in \Z$, by the same computation as above, we have
that $h(T^nx_0)=e(nt)f(\tau^n\cdot x_1)$ and we are done.
\end{proof}

\begin{corollary}
\label{cor:NorthoAP}
Let $(X,T)$ be an elementary nilsystem. Then 
$\Nil(X)\perp\AP$, meaning that for every $\ua\in\Nil(X)$ and every 
$\ub\in\AP$, we have $\Av(\ua\,\ub)=0$.

If $\ua$ is an elementary nilsequence and is not almost periodic, 
then $\ua\perp\AP$.
\end{corollary}

\begin{proof}
Every function in $\CC_1(X)$ obviously has zero integral with 
respect to the Haar measure of $X$.  Thus by definition and unique 
ergodicity, every sequence in $\Nil(X)$ has zero average.
The first statement follows immediately from Part~\eqref{it:invarNilX} of 
Lemma~\ref{lem:basepoint}. The second statement follows by density.
\end{proof}

\subsection{Construction of the Bohr extension}
\label{subsec:constr-universal}
Let $(X=G/\Gamma,T)$ be a minimal elementary $2$-step nilsystem 
in reduced form 
$X=G/\Gamma$.  We use following notation.

$T$ denotes the translation by $\tau\in G$, 
$\pi\colon X\to Z=G/G_2\Gamma$ and $p\colon G\to Z$ are the natural 
projections, $\sigma=p(\tau)$ and $S$ is the translation by $\sigma$ 
on $Z$.  (So $\pi\colon (X,T)\to(Z,S)$ is a factor map.)
 $Z$ is endowed with its Haar measure $m_Z$. 
Let $\BZ$ denote the Bohr compactification of $\Z$: it is the dual of the 
circle group endowed with the 
discrete topology and it contains $\Z$ as a 
dense subgroup.  Let $\BZ$ be endowed with its Haar measure $m_{\BZ}$.
We write $R$ for the (minimal) translation by $1$ in 
$\BZ$.  There exists a (well defined) continuous group homomorphism 
$r\colon\BZ\to Z$ such that $r(1)=\sigma$. This homomorphism is 
onto, and it is a factor map from $(\BZ,R)$ to $(Z,S)$.

Define
$$
 \tilde X=\{(x,z)\in X\times \BZ\colon \pi(x)=r(z)\}
$$
and let $\tilde T$ be the restriction of $T\times R$ to $\tilde X$. 
We write $q_1\colon\tilde X\to X$ and $q_2\colon\tilde X\to\BZ$ for 
the natural projections.

The system $(\tilde X,\tilde T)$ is called the \emph{Bohr extension} 
of $(X,T)$.

We give a second presentation of this system.
Define:
$$
 \tilde G=\bigl\{(g,z)\in G\times\BZ\colon p(g)=r(z)\bigr\}\ ;\ 
\tilde\Gamma=\bigl\{ (\gamma,0)\colon \gamma\in\Gamma\}
$$
and define $\tilde\tau=(\tau,1)\in G\times\BZ$.
Then $\tilde G$ is a closed subgroup of $G\times\BZ$, 
$\tilde\tau\in\tilde G$ and $\tilde\Gamma$ is a discrete subgroup 
of $\tilde G$. The map $(g,z)\mapsto (p(g),z)$ from $G\times\BZ$ to
$X\times\BZ$ induces a homeomorphism of $\tilde G/\tilde\Gamma$ 
onto $\tilde X$ and we identify these spaces. Under this 
identification, $\tilde T$ is the translation by $\tilde\tau$ on 
$\tilde X$.

We have that $\tilde G_2=\{(u,0)\colon u\in G_2\}$ and thus can be 
identified with the circle group $\S^1$. The action of this group on 
$\tilde X$ is given by
$$
 u\cdot(x,z)=(u\cdot x,z)\quad(u\in \S^1=G_2,\ (x,z)\in\tilde X)\ .
$$
This action of $G_2$ commutes with $\tilde T$ and 
 the quotient of $\tilde X$ by this action can be
identified with $\BZ$ in a natural way, the identification being 
given by the projection $q_2$.

The next result follows classically from the fact that all  eigenfunctions 
of $X$ factorize through $Z$ and that $(\bar Z,R)$ is a rotation. 
For completeness, we give a proof in a appendix~\ref{appendix:A}.

\begin{lemma}
\label{lem:tildeX}
$(\tilde X,\tilde T)$ is uniquely ergodic and the topological support 
of its invariant measure is equal to $\tilde X$.
\end{lemma}

\begin{corollary}
$(\tilde X,\tilde T)$ is minimal.
\end{corollary}

\subsection{Applications to elementary nilsequences}
\label{subsec:uses-for-elementary}

We keep the notation of the previous sections.
As for elementary nilsystems, we define
$\CC_1(\tilde X)$ to be the space of continuous functions $f$ on $\tilde 
X$ such that
$$
f(u\cdot x,w)=uf(x,w)\quad (u\in \S^1,\ (x,w)\in\tilde X)\ .
$$

\begin{lemma}
\label{lem:nilXtilde} 
Let $\tilde x_0\in\tilde X$.
The map 
\begin{equation}
\label{eq:nilXtilde}
 f\mapsto\bigl(f(\tilde T^n\tilde x_0\colon n\in\Z\bigr)
\end{equation}
is an isometry of $\CC_1(\tilde X)$ onto the space $\Nil(X)$.
\end{lemma}
\begin{proof}
Let $x_0=q_1(\tilde x_0)$. We identify $\wh{\BZ}$ with $\T$ (with the 
discrete topology).

Clearly, $\CC_1(\tilde X)$ is the closed 
(under the uniform norm) linear span of the family of functions 
$(x,w)\mapsto \phi(x)\chi(w)$ for $\phi\in\CC_1(X)$ and 
$\chi\in\wh{\BZ}$.  If $f$ is a function of this type, then for every 
$n\in\Z$, we have
$$
 f(\tilde T^n\tilde x_0)=\phi(T^nx_0)e(n\chi)\ .
$$
Thus the sequence $( f(\tilde T^n\tilde x_0))$ belongs to 
$\Nil(X)$ by part~\eqref{it:invarNilX} of Lemma~\ref{lem:basepoint}.
By density, the same property remains valid for every 
$f\in\CC_1(\tilde X)$ and thus formula~\eqref{eq:nilXtilde} 
defines a map $\jmath\colon \CC_1(\tilde X)\to\Nil(X)$.

If $f\in\CC_1(\tilde X)$ and $\ua=\jmath(f)$, by minimality of 
$(\tilde X,\tilde T)$ we have that
$$
 \norm{\ua}_\infty=\sup_{n\in\Z}| f(\tilde T^n\tilde x_0)|=\norm{f}_\infty
$$
and the map $\jmath$ is an isometry. 

By Part~\eqref{it:basepoint} of Lemma~\ref{lem:basepoint},
$\Nil(X)$ is the closed linear span 
of the family of sequences of the type $(e(nt)h(T^nx_0))$, where 
$t\in\T$ and $h\in\CC_1(X)$.  We are left with showing that every
sequence of this type belongs to the range of $\jmath$.
Let $\chi$ be the character of $\BZ$ 
corresponding to $t$ and let $f(x,w)=h(x)\chi(w)$. This function belongs 
to $\CC_1(\tilde X)$, and its image under the isometry $\jmath$ is 
the given sequence.
\end{proof}

\begin{corollary}
\label{cor:NXirreducible}
Let $(X,T)$ be an elementary nilsystem. Then $\Nil(X)$ is irreducible, 
meaning that it does not contain any closed proper nontrivial 
subspace that is invariant under the shift and multiplication by linear 
exponential sequences.

In particular, for every non-zero $\ua\in\Nil(X)$, we have that 
$\CS(\ua)=\Nil(X)$.
\end{corollary}

\begin{proof}
Let $\ua$ a non-identically zero elementary nilsequence 
belonging to $\Nil(X)$. We need to show that the class $\CS(\ua)$ 
of this sequence is equal $\Nil(X)$.

Through the isometry defined in Lemma~\ref{lem:nilXtilde}, 
the sequence $\ua$ corresponds to 
a non-identically zero function $f$ on $\tilde X$ that belongs to 
$\CC_1(\tilde X)$. 
Since for every $t\in\T$ and every $k\in\Z$ the sequence 
$\ue(t)\sigma^k\ua$ corresponds to the function $\chi\circ p_2\cdot 
T^kf$, where $\chi$ is the character of $\BZ$ associated to $t$,
the space  $\CS(\ua)$ corresponds through this isometry 
to the closed linear 
subspace $\CF$ of $\CC_1(\tilde X)$ spanned by the functions 
$\chi\circ q_2\cdot f\circ T^n$, where $\chi\in\wh{\BZ}$ and 
$n\in\Z$.  It remains to show that this space is equal to $\CC_1(X)$.

Let $U$ be the non-empty open subset 
$\{x\in\tilde X\colon f(x)\neq 0\}$ of $\tilde X$.
As $f\in\CC_1(\tilde X)$ and the quotient of $\tilde X$ 
under the action of $\S^1$ is equal to $\BZ$, we have that 
$U=p_2\inv(V)$ where $V$ is some open subset of $\BZ$.

Since the closed linear span of the functions 
$\chi\circ p_2$ for $\chi\in\wh{\BZ}$  is equal to 
$\{h\circ p_2\colon h\in\CC(\BZ)\}$, we have that 
 the closed linear span $\CF$ of the family $\chi\circ q_2\cdot f$ 
for $\chi\in\wh{\BZ}$ is equal to 
$$\{g\in\CC_1(\tilde X)\colon g=0\text{ outside of }U\}\ .$$

By minimality, there exists an integer $m\geq 1$ such that
the sets $R^jV$, $0\leq j\leq m$, cover $\BZ$. 
We chose a partition of the unity $\{\phi_j\colon 0\leq j\leq m\}$ in
 $\CC(\BZ)$, subordinated to this cover. 
For any $h\in\CC_1(\BZ)$, each of the functions $h\cdot\phi_j\circ 
p_2$ belongs to $\CF$, and thus $h$ belongs to $\CF$.
\end{proof}

In this proof, the shift was only needed because the 
function $f$ on $\tilde X$ may vanish at some point.

\begin{corollary}
\label{cor:factorelem}
Let $(X=G/\Gamma,T)$ and $(X'=G'/\Gamma',T')$ be two elementary nilsystems and assume 
that there exists a factor map
$q\colon (X,T)\to (X',T')$  commuting with the action of 
$\S^1=G_2=G'_2$ on these systems. Then $\Nil(X)=\Nil(X')$.
\end{corollary}

\begin{proof}
The factor map induces an inclusion of $\Nil(X)$ in $\Nil(X')$ and 
equality follows from the irreducibility of $\Nil(X')$.
\end{proof}

\begin{remark*}
One can check that in this case $\widetilde{X'}=\tilde X$ and this 
leads to another proof of the same result.
\end{remark*}

\begin{corollary}
\label{cor:SaNx}
Let $\ua$ be a non-zero basic elementary nilsequence.  Then
there exists an elementary nilsystem $X$ such that 
$\CS(\ua)=\Nil(X)$.
\end{corollary}

\section{Taking products}
\label{sec:products}
We know that the product of two nilsequences is a nilsequence. We 
consider here the case of two elementary nilsequences.
If $\CA$ and $\CB$ are two families of sequences, we write
$$
 \CA\cdot\CB:=\{\ua\,\ub\colon\ua\in\CA,\ \ub\in\CB\}\ .
$$

Let $(X,T)$ be an elementary nilsystem. 
Then  for every $\phi,\psi\in\CC_1(X)$, the function $\phi\bar\psi$ 
belongs to $\CC_0(X)$ and for every  every $x_0\in X$, the sequence 
$(\phi(T^nx_0)\overline{\psi(T^nx_0)}\colon n\in\Z)$ 
is almost periodic. Thus
$$
 \Nil(X)\cdot\overline{\Nil(X)}\subset\AP\ .
$$
\subsection{Orthogonality}
\label{subsec:orthogonality}

\begin{proposition}
\label{prop:equalortho2}
 If $(X=G/\Gamma,T)$ and $(X'=G'/\Gamma',T')$ are two elementary 
nilsystems, then exactly one of the following statements holds:
\begin{itemize}
\item $\Nil(X) = \Nil(X')$, and in this case $\Nil(X)\cdot\overline{\Nil(X')} 
\subset \AP$.
\item $\Nil(X) \perp \Nil(X')$, 
and in this case there exists and elementary nilsystem $(Y,S)$ 
such that $\Nil(X)\cdot\overline{\Nil(X')} \subset N(Y)$.  
\end{itemize}
More precisely, $\Nil(X)=\Nil(X')$ if and only if there exists a closed subgroup 
$H$ of $G\times G'$ such that
\begin{enumerate}
\item
  \label{it:tautauprme}
$(\tau,\tau')\in H$.
  \item
  \label{it:diagonal}
The  commutator subgroup $H_2$ of $H$ is the diagonal subgroup of
$\T^2=G_2\times G'_2$.
  \item
 \label{it:cocompactLambda}
The subgroup $\Lambda=(\Gamma\times\Gamma')\cap H$
 of $H$ is cocompact in $H$.
\end{enumerate}
\end{proposition}

The invariance under products of the family of elementary nilsequences 
follows immediately. Indeed, the product of an elementary nilsequence 
by an almost periodic sequence is a nilsequence; 
Proposition~\ref{prop:equalortho2} shows that the product of two 
basic elementary nilsequences is an elementary nilsequence, and the
general case follows by density.

\begin{proof}
\strut
\subsubsection{} 
\label{subsec:notexist}
First we assume that there exists no group $H$ satisfying
properties~\eqref{it:tautauprme},  \eqref{it:diagonal} 
and~\eqref{it:cocompactLambda} 
above.

Let $x_0\in X$ be the image  if the unit 
element of $G$ and $x'_0\in X'$ the image of the unit element of $G'$.
Then $\Gamma\times\Gamma'$ is the stabilizer
of $(x_0,x'_0)$ under the action of $G\times G'$ on $X\times 
X'$.

Let $W$ be the closed orbit under $T\times T'$ of the point 
$(x_0,x'_0)$ in the nilsystem $(X\times X', T\times T')$. 
It is classical that $(W,T\times T')$ is a two step nilsystem. 
Write $W=H/\Lambda$, where $H,\Lambda$ satisfy 
properties~\eqref{it:tautauprme} and~\eqref{it:cocompactLambda} above.
By hypothesis, property~\eqref{it:diagonal} is not satisfied and 
$H_2$ is not equal to the diagonal subgroup of $\T^2$.

By minimality, the natural projection $W\to X$ is onto. 
It follows that the image of the natural projection  $H\to G$ 
has countable index in $G$. 
Therefore this image is open in $G$ and thus contains $G_0$.  
As it contains $\tau$, it is equal to $G$. 
For the same reason, the natural projection $H\to G'$ is 
onto.

Therefore $H$ is not abelian. Its commutator subgroup is therefore 
not the trivial subgroup of $\T^2$. 
Let $\chi$ be the restriction to $H_2$ of the character 
$(u,u')\mapsto u-u'$ of $\T^2$.  Thus $\chi$ is not the trivial 
character of $H_2$.

Let $(Y,S)$ be the elementary nilsystem obtained by taking the quotient of  $X$ 
by the subgroup $\ker(\chi)$ of $G_2$ as in 
Section~\ref{subsec:decomposition}.  Let $p\colon W\to Y$ be
the natural projection and set $y_0=p(x_0,x'_0)$.

Let $\phi\in\CC_1(X)$ and $\phi'\in\CC_1(X')$. 
Let $\Phi$ be the restriction to $Y$ of the function 
$(x,x')\mapsto\phi(x)\overline{\phi'(x')}$.
This function belongs to $\CC_\chi(W)$ and thus can be written 
as $\Psi\circ p$ for some function $\Psi\in\CC_1(Y)$. 
For every $n\in\Z$, we have
$$
\phi(T^nx_0)\overline{\phi'((T')^nx'_0)}=\Phi((T\times T')^n(x_0,x'_0))
=\Psi(S^ny_0)\ .
$$
Since $\Nil(X)$ is the closed linear span of the family of sequences of 
the form $(\phi(T^nx_0)e(nt))$, and similarly for $\Nil(X')$, we have that
$\Nil(X)\cdot\overline{\Nil(X')}\subset\Nil(Y)$. 
By irreducibility we have that 
$$
 \text{the closed linear span of }\Nil(X)\cdot\overline{\Nil(X')}
\text{ is equal to }\Nil(Y)\ .
$$
Since every sequence in $\Nil(Y)$ has zero average,
 $\Nil(X)\perp\Nil(X')$.

\subsubsection{} We now assume that there exists a subgroup $H$ of $G\times 
G'$ with properties~\eqref{it:tautauprme}, \eqref{it:diagonal} 
and~\eqref{it:cocompactLambda}.

 Let $L$ be the closed subgroup of $G\times G'$ spanned by $\Lambda$
and $\sigma=(\tau,\tau')$.

Set $W=L/\Lambda$ and let $S\colon W\to W$ be the translation by 
$\sigma=(\tau,\tau')$. 
Then $(W,S)$ is a $2$-step nilsystem and the 
natural projections $W\to X$ and $W\to X'$ are factor maps. 
By the same argument as in Section~\ref{subsec:notexist}, the natural 
projections $L\to G$ and $L\to G'$ are onto and $L$ is not abelian.
The commutator subgroup of $L$ is thus a nontrivial subgroup of 
$H_2$, and by hypothesis, it is the diagonal subgroup of $\T^2$.
Since $L_2$ is connected, $L_2=H_2$.

Substituting $L$ for $H$, we are reduced to the case that
\begin{enumerate}
\setcounter{enumi}{3}
\item
\label{it:HLambda}
$H$ is the closed subgroup of $G\times G'$ spanned by $\Lambda$ and
$\sigma$.
\end{enumerate}

 It follows that $(W,S)$ is an elementary nilsystem. 
The natural factor maps $W\to X$ and $W\to X'$ commute with the 
actions of $\T=G_2=G'_2=H_2$ on $X,X'$ and $W$, respectively.
By Corollary~\ref{cor:factorelem}, we have that 
$\Nil(X)=\Nil(W)=\Nil(X')$.
\end{proof}

\subsection{} For later use, we show:
\begin{lemma}
\label{lem:Phi}
Let $X,X'$ be as in Proposition~\ref{prop:equalortho2} and let $H$  be a 
closed subgroup of $G\times G'$ satisfying the 
conditions~\eqref{it:tautauprme},
~\eqref{it:diagonal} and~\eqref{it:cocompactLambda} of 
Proposition~\ref{prop:equalortho2}.  
Then there exist a subgroup $\Gamma_1$ of $\Gamma$ with finite index, a 
subgroup $\Gamma'_1$ of $\Gamma'$ with finite index and an isomorphism
$\Phi\colon\Gamma_1\to\Gamma'_1$ such that
\begin{gather}
\label{eq:Phicommut}
 \text{for every }\gamma\in\Gamma_1,\quad 
[\Phi(\gamma),\tau']=[\gamma,\tau]\ ;\\
\label{eq:LambdaPhi}
\Lambda=\bigl\{\bigl(\gamma,\Phi(\gamma)\bigr)\colon\gamma\in\Gamma_1\bigr\}\ .
\end{gather}
\end{lemma}

\begin{proof}
Let $\Gamma_1$ and $\Gamma'_1$ be the images of $\Lambda$ under the 
natural projections $H\to G$ and $H\to G'$, respectively.

By~\eqref{it:cocompactLambda}, $\Lambda$ is cocompact in $H$. We have 
shown above that the natural projection $H\to G$ is onto. It follows 
that  $\Gamma_1$ is cocompact in $G$ and thus it has finite index in 
$\Gamma$. By the same proof, $\Gamma'_1$ has finite index in 
$\Gamma'$.

For every $\gamma\in \Gamma_1$, by definition 
there exists $\gamma'\in\Gamma'_1$ with $(\gamma,\gamma')\in\Lambda$. 
Properties~\eqref{it:tautauprme} and~\eqref{it:diagonal} 
imply that 
$[\gamma',\tau']=[\gamma,\tau]$. Since the map $\gamma'\mapsto 
[\gamma',\tau']$ is one to one, this condition completely determines 
$\gamma'$ for a given $\gamma$.  Thus there exists a group 
homomorphism $\Phi\colon\Gamma_1\to\Gamma'$ 
satisfying~\eqref{eq:Phicommut} and~\eqref{eq:LambdaPhi}.
By exchanging the role played by the two coordinates, we have that 
$\Phi$ is a group isomorphism from $\Gamma_1$ onto $\Gamma'_1$.
\end{proof}

\subsection{Proof of Theorem~\ref{th:equalortho}}
\label{subsec:proofequalortho}
We start with Lemma.
\begin{lemma}
\label{lem:basicelement}
For every non-almost periodic elementary nilsequence $\ua$, 
there exist a basic 
elementary nilsequence $\ub$ with $\CS(\ua)=\CS(\ub)$ and an 
elementary nilsystem $(X,T)$ with $\CS(\ua)=\Nil(X)$.
\end{lemma}

\begin{proof}
Let $\ua$ be a non-almost periodic  elementary nilsequence; in 
particular, $\ua$ is not identically zero.
There exists a sequence $(\ub(j)\colon j\geq 1)$ of basic elementary nilsequences 
converging uniformly to $\ua$. For every $j\geq 1$, let $(X_j,T_j)$ be 
an elementary nilsystem such that $\ub(j)$ belongs to $\Nil(X_j)$.

As the quadratic norm  is a norm 
(Section~\ref{subsec:quadratic-norm}),
 $\norm{\ua}_2>0$ and 
$\norm{\ua-\ub(j)}_2\to 0$ as $j\to+\infty$.
Therefore, there exists $i\geq 1$ such that 
$\norm{\ub(i)-\ub(j)}_2<\norm{\ub(i)}_2$ for every $j>i$.

This implies that for $j>i$ the sequence $\ub(j)$ is not 
orthogonal to the sequence $\ub(i)$.
The spaces  
$\Nil(X_j)$ and $\Nil(X_i)$ are not orthogonal and by
Proposition~\ref{prop:equalortho2}, these spaces are equal.
Therefore $\ua(j)$ belongs to $\Nil(X_i)$ for every $j>i$ and 
so $\ua\in\Nil(X_i)$.  By Corollary~\ref{cor:SaNx}, 
this space is equal to $\CS(\ub(j))$ 
\end{proof}

We now prove Theorem~\ref{th:equalortho} by collecting the 
results of the preceding Sections.

\eqref{it:equalortho}
Let $\ua$ and $\ub$ be two non-identically zero elementary 
nilsequences. We want to show that the spaces $\CS(\ua)$ and $\CS(\ub)$ 
are equal or orthogonal. If one of these sequences, say $\ua$, is 
almost periodic then $\CS(\ua)=\AP$
and the result follows from Corollary~\ref{cor:NorthoAP}.

Assume now that $\ua$ and $\ub$ are not almost periodic.
By Lemma~\ref{lem:basicelement} and 
Corollary~\ref{cor:NXirreducible},
 there exist two 
elementary nilsystems $(X,T)$ and $(X',T')$ with $\CS(\ua)=\Nil(X)$ 
and $\CS(\ub)=\Nil(X')$ and the result follows from
Proposition~\ref{prop:equalortho2}.

\eqref{it:open} The same proof used for Lemma~\ref{lem:basicelement} 
shows that every class is open. 

\eqref{it:products} Follows directly from Lemma~\ref{lem:basicelement} and 
Proposition~\ref{prop:equalortho2}.

\eqref{it:groups} The unit element for the multiplication is $\AP$. 
Each class different from $\AP$ is equal to $\CN(X)$ for some elementary 
nilsystem $X$, and by Proposition~\ref{prop:equalortho2},
the inverse of this class is $\CN(\bar X)$.\qed

\section{Classification: reduction to the connected case}
Here we begin the proof of Theorem~\ref{th:classification}, 
considering first the case of an elementary system arising from a 
non-connected group.
\begin{proposition}
\label{prop:non-connected}
Let $(X=G/\Gamma,T)$ be a (minimal) elementary nilsystem. Then there exist 
$t\in\T$ and an elementary nilsystem $(X'=G'/\Gamma',T')$ with $G'$ 
connected such that $\Nil(X)=\uq(t)\cdot\Nil(X')$.
\end{proposition}

\begin{proof}
Assume that $(X = G/\Gamma, T)$ is an elementary nilsystem 
and assume that $T$, as usual, is rotation by the element $\tau\in G$.  
Let $G_0$ denote the connected component of the 
identity in $G$ and set $\Gamma_0 = \Gamma\cap G_0$ and 
$X_0 = G_0\Gamma/\Gamma_0\cong G_0/\Gamma_0$.  
We have that $X_0$ is a connected component of $X$ and is open in $X$.  

\subsection*{First step}
We first build a particular element 
 $\tau_0\in G_0\Gamma$ 
such that the system $(X_0, T_0)$, where $T_0$ is translation 
by $\tau_0$, is minimal.  

Since $X$ is minimal, there exists $k\in\N$ such that $T^kX_0\cap X_0$ 
is nonempty.  Choose $k$ to the be least integer such that this 
intersection is nonempty.  
(If $X$ is connected, then $X=X_0$ and $k=1$.)
By definition of $X_0$, $K$ is also the smallest integer with 
$\tau^k\in G_0\Gamma$ and $T^kX_0 = X_0$.  
Since $G$ is spanned by $G_0$ and $\tau$, it follows that
$G_0\Gamma$ is spanned by $G_0$ and $\tau^k$. 
Since $(X,T)$ is minimal, the set
$$
\{\tau^n\gamma\colon n\in\Z, \gamma\in\Gamma\}\cap G_0\Gamma = 
\{\tau^{kn}\colon n\in\Z, \gamma\in\Gamma\}
$$
is dense in $G_0\Gamma$.  
It follows that $(X_0, T^k)$ is a minimal system.  

Since $\tau^k\in G_0\Gamma$ and $G_0$ is connected, there 
exists $\theta\in\Gamma$ and $\tau_1\in G_0$ such that 
\begin{equation}
\label{eq:theta-k}
\tau^k = \tau_1^k\theta \ .
\end{equation}

We now consider two cases, depending on whether $G_0$ is abelian or not.

\subsubsection*{The abelian case}

First assume that $G_0$ is abelian.  
Then $X_0$ is a compact connected abelian 
group of finite dimension and so is a torus.  Furthermore, 
$X_0$ contains the one dimensional torus $G_2 = \S^1$ as 
a closed subgroup, and so $X_0$ can be identified with 
$Z_0\times G_2$, where 
$$
Z_0 = G_0\Gamma/G_2\Gamma\cong G_0/G_2\Gamma_0\ .
$$
Here $Z_0$ is an open and connected subgroup of 
$Z = G/G_2\Gamma$ and so $Z_0$ is isomorphic to $\S^d$ 
for some $d\in\N$.  
Let $\pi\colon X\to Z$ denote the natural projection and 
so $\pi(\tau)$ denotes the image of $\tau$ in $Z$.  
Let $\alpha$ denote the image of the element $\tau_1\in G_0\Gamma$ 
in $Z_0$.  By Equation~\eqref{eq:theta-k}, $\alpha^k = \pi(\tau)^k$.  
Since $(X_0, T^k)$ is minimal, the translation by 
$\alpha^k$ on $Z_0$ is minimal.  
Therefore, there exists $v\in G_2$ such that the translation 
by $(\alpha^k, v)$ on $Z_0\times G_2$ is minimal.
Since $G_0$ is connected, there exists $u\in G_2$ such that 
$v = u^k$.  

Let $\tau_0$ be a lift of $(\alpha, u)$ in $G_0$ and 
let $T_0$ be the translation by $\tau_0$ on $X_0$.  
Then $(X_0, T_0)$ is a rotation system and $(X_0, T_0^k)$ 
is minimal.  It follows that $(X_0, T_0)$ is minimal.  
Since $\tau_0^k$ and $\tau^k$ both project to $\pi(\tau)^k$ 
on $Z$, there exists $\eta\in\Gamma$ and $w\in G_2$ such that 
$\tau^k = \tau_0^k\eta w$.  Thus 
we have constructed an element $\tau_0$ of $G_0$ such that:
\begin{enumerate}
\item 
\label{enum:one}
$\tau^k = \tau_0^k\eta w$ for some $\eta\in\Gamma$ and $w\in G_2$;
\item The system $(X_0, T_0)$, where $T_0$ is  the translation by $\tau_0$, 
is minimal;
\item 
\label{enum:three}
The system $(X_0, T_0)$ is either a rotation on a compact 
abelian group or is an elementary nilsystem.
\end{enumerate}

\subsubsection*{The non abelian case}
We now turn to the case that $G_0$ is not abelian.
We take $\tau_0= \tau_1$, 
$\eta = 1$ and $w = 1$ and show that the 
properties~\eqref{enum:one}--\eqref{enum:three} are satisfied.

Since $G_0$ is now assumed to be nonabelian, 
we have that $(G_0)_2$ is a nontrivial torus.  
It is a subgroup of $G_2=\S^1$ and $(G_0)_2 = G_2$.  

Again, $G$ is spanned by $G_0$ and $\tau$.  Let $\gamma\in\Gamma$ 
and write $\gamma = g_0\tau^j$ for some $g_0\in G_0$ and $j\in\Z$. 
By definition, $\Gamma$ is the stabilizer of some point of $X$ 
and so $j = km$ for some $m\in\Z$.   Thus
$$
\gamma = g_0\tau^{km} = 
g_0\tau_0^{km}[\tau_0, \theta]^{km(m-1)/2}\theta^m 
= g_1\theta^m
$$
for some $g_1\in G_0$.  
Therefore $g_1\in G_0\cap\Gamma = \Gamma_0$ and we conclude that 
$\Gamma$ is spanned by $\Gamma_0$ and $\theta$.  
Recall that the set $\{\tau^{kn}\gamma\colon n\in\Z, \gamma\in\Gamma\}$ 
is dense in $G_0\Gamma$.  By Equation~\eqref{eq:theta-k}, 
this set is included in 
$$
\{\tau_0^{kn}u\gamma\colon n\in \Z, u\in G_2, \gamma\in\Gamma\}
$$
and so this set is also dense in $G_0\Gamma$.  
Since $G_0$ is open, the set 
$$
\{\tau_0^{kn}u\gamma_0\colon n\in \Z, u\in G_2, \gamma\in\Gamma\}\cap G_0 
 = 
\{\tau_0^{kn}u\gamma_0\colon n\in \Z, u\in G_2, \gamma\in\Gamma_0\}
$$ 
is dense in $G_0$.  
Since $G_0$ is connected and $(G_0)_2 = G_0$, we 
have that $(X_0 = G_0/\Gamma_0, T_0^k)$ is minimal.  
Thus $(X_0, T_0)$ is an elementary nilsystem 
satisfying properties~\eqref{enum:one}--\eqref{enum:three}.

\subsection*{Second step}

Finally we compare the spaces $\Nil(X)$ and $\Nil(X_0)$.
Let $x_0$ be the image in $X$ of the unit element of $G$.
The space $\Nil(X_0)$ is spanned by translates of sequences of 
the form 
$$
\ua = (\phi(T^nx_0)\colon n\in\Z) \ ,
$$
where $\phi\in \C_1(X)$ vanishes outside of $X_0$.  
By definition of 
$T_0$, we have $a_n = 0$ if $n\notin k\Z$.  
We also use $\phi$ to denote the restriction of $\phi$ to $X_0$ 
and define the sequence $\ub$ by 
$$
b_n =  \frac{1}{k}\sum_{j=0}^{k-1}e(jn)\phi(T_0^nx_0)
=  \begin{cases}
\phi(T_0^nx_0) &\text{ if }n\in k\Z\\
0 &\text{ otherwise }\ .
\end{cases}
$$

By construction, $\tau^k = \tau_0^k\eta w$ for some $\eta\in\Gamma$
and $w\in G_2$.  Let $u\in G_2$ be defined by 
$\theta\cdot x_0 = u\cdot x_0$.  
If $n = km$ for some $m\in\Z$, we have that 
\begin{multline*}
a_n = \phi((\tau_0^k\theta w)^m\cdot x_0) = 
\phi(\tau_0^{km}[\tau_0,\eta]^{km(m-1)/2}w^m\theta^m \cdot x_0)
\\
= 
[\tau_0, \eta]^{km(m-1)/2}\,w^mu^n\,\phi(\tau_0^{km}\cdot x_0) 
= e_n(s)q_n(t)b_n 
\end{multline*}
for some $s, t\in\T$ which do not depend on choice of $\phi$.  

For $n\notin k\Z$, we have that $a_n = 0 = e_n(s) q_n(t)b_n$ 
and so $\ua = \uq(t)\ue(s)\ub$.  Thus
$\ua\in \uq(t)\Nil(X_0)$.  Since 
$\Nil(X)$ and $\uq(t)\Nil(X_0)$ are irreducible, 
we have that $\Nil(X) = \uq(t)\Nil(X_0)$.  
In conclusion, we have constructed an elementary nilsystem such that 
$\Nil(X) = \uq(t)\Nil(X')$.  
\end{proof}

\section{Classification: The connected case}

\subsection{Reduction to Heisenberg systems}
For Heisenberg systems we use the definitions and notation  
of Section~\ref{subsec:defHeisenberg}. The condition for minimality  
of Heisenberg nilsystems was given in Section~\ref{subsec:backexamples}.

Throughout this section, we consider $\Nil(X)$ for various spaces and 
various transformations on $X$.  To minimize confusion, 
we write $\Nil(X,T)$ instead of the usual $\Nil(X)$.  

\begin{lemma}
\label{lem:connected-Heisen}
Let $(X=G/\Gamma,T)$ be an elementary nilsystem (minimal, written in 
reduced form) where $G$ is connected.
Then there exists an Heisenberg minimal system $(N_d,T')$ 
with $\Nil(X,T)=\Nil(H_d,T')$.
\end{lemma}

\begin{proof}
By Lemma~\ref{lem:desc-connected} in Appendix~\ref{sec:desc-conected}, 
$(X,T)$ is the product of of an elementary system $(X', T'')$ 
of the form given by the lemma and a rotation.  Since
$\Nil(X')$ is invariant under taking products with almost periodic sequences, 
it is easy to check that $\Nil(X,T)=\Nil(X',T'')$.

We can therefore reduce to the case that 
$G=\R^{2d}\times \S^1$ and $\Gamma=\Z^{2d}\times \{1\}$
for some $d\geq 1$, with multiplication 
given by 
$$
 (x,z)\cdot(x',z')=(x+x',zz'e(\langle Ax\mid x'\rangle)\quad
(x,x'\in\R^{2d},\ z,z'\in\S^1)
$$
where $A$ is a $2d\times 2d$ matrix such that the matrix $B=A-A^t$ is 
nonsingular. We build a Heisenberg nilsystem $(N_d,T')$ with 
$\Nil(X,T)=\Nil(N_d,T')$.

Let $J_{2d}$ be as in Section~\ref{subsec:classific}.
Choose a $2d\times 2d$ matrix $\Phi$ with rational 
entries such that 
$$
 \Phi^tJ_{2d}\Phi=B\ .
$$
(Existence of such a matrix follows using an antisymmetric 
version of Gauss decomposition of a quadratic form.) 
Let $\tau=(\delta,\gamma)\in\R^{2d}\times\S^1=G$ be the element 
defining $T$, $\tau'=(\Phi(\delta),1)\in \R^{2d}\times\S^1=H_d$ and 
let $T'$ be the translation by $\tau'$ on $N_d = 
H_d/(\Z_d\times\Z_d\times\{1\})$.

Since $(X,T)$ is minimal, the coordinates of $\delta$ are rationally 
independent modulo $1$ by Proposition~\ref{prop:cond-min}.  Since
$\Phi$ has rational entries, the coordinates of $\Phi(\delta)$ are 
rationally independent modulo $1$ and $(N_d,T')$ is minimal. 
Define 
$$
 H=\bigl\{ \bigl((x,z),(x',z')\bigr)\in G\times H_d\colon 
 x'=\Phi(x),\ z'=z\bigr\}\ .
$$
Then $H$ is a subgroup of $G\times N_d$ satisfying the conditions of 
Proposition~\ref{prop:equalortho2} and thus $\Nil(X,T)=\Nil(N_d,T')$.
\end{proof}

\subsection{The case of Heisenberg systems}

\begin{lemma}
\label{lem:class-Heisen}
Let $(N_d=H_d/\Lambda_d,T)$ and $(N_{d'}= H_{d'}/\Lambda_{d'},T')$ 
be two minimal Heisenberg systems, 
where the transformations $T$ and $T'$ are respectively the 
translation by 
$\tau=(\alpha,\beta,\gamma)\in\R^d\times\R^d\times\S^1=H_d$ and 
translation by $\tau'=(\alpha',\beta',\gamma')
\in \R^{d'}\times\R^{d'}\times\S^1=H_{d'}$.

Then $\Nil(N_d,T)=\Nil(N_{d'},T')$ if and only if $d=d'$ and there 
exists a matrix $Q\in\Sp_{2d}(\Q)$ mapping 
$(\alpha,\beta)$ to 
$(\alpha',\beta')$.
\end{lemma}
(We consider $Q$ as a linear map from $\R^d\times\R^d$ to itself.)
\begin{proof}
First assume that $\Nil(N_d,T)=\Nil(N_{d'},T')$. Let $H,\Lambda$ be 
as in Proposition~\ref{prop:equalortho2} and $\Gamma_1,\Gamma'_1,\Phi$ 
as in Lemma~\ref{lem:Phi}. Since $\Gamma_1$ is of finite index in the 
group $\Lambda_d$ which is isomorphic to $\Z^d$, it is itself 
isomorphic to $\Z^d$. 
In the same way, $\Gamma'_1$ is isomorphic to 
$\Z^{d'}$. Since $\Phi\colon\Gamma_1\to\Gamma'_1$ is an isomorphism, $d=d'$.
Moreover, $\Phi$ is given by a $2d\times 2d$ matrix $Q$ with rational 
entries. Relation~\eqref{eq:Phicommut} of Lemma~\ref{lem:Phi} 
means exactly that $Q$ belongs to $\Sp_{2d}(\Q)$, and 
relation~\eqref{it:tautauprme} of Proposition~\ref{prop:equalortho2}
that $Q$ maps 
$(\alpha,\beta)$ to 
$(\alpha',\beta')$.

Conversely, assume that $d=d'$ and that the matrix $Q\in\Sp_{2d}(\Q)$ 
maps $(\alpha,\beta)$ to 
$(\alpha',\beta')$.
 We set 
 $$
 H=\bigl\{\bigl((x,y,z),
(x',y',z')\bigr)\in N_d\times N_d\colon
(x',y')=Q(x,y),\ z'=z\bigr\}
$$
where $x,y,x',y'\in\R^d$ and $z,z'\in\R$. Then $H$ satisfies all the conditions of 
Proposition~\ref{prop:equalortho2} and thus $\Nil(N_d,T)=\Nil(N_{d'},T')$.
\end{proof}

\subsection{Heisenberg systems and quadratic exponential sequences}

\begin{lemma}
\label{lem:convsnot}
Let $(N_d,T)$  be a minimal Heisenberg nilsystem and 
let $t\in\T$ be such that for any integer $m>0$, $mt$ 
does not belong to the eigenvalue group of $(N_d,T)$. Then 
the sequence $\uq(t)$ does not belong to $\Nil(N_d,T)$  
and there exists an elementary nilsystem $(X,S)$ with
$\uq(t)\cdot \Nil(N_d,T)=\Nil(X,S)$.
 
Moreover, for every Heisenberg system $(N_{d'},T')$, 
$\uq(t)\cdot\Nil(N_d,T)\neq \Nil(N_{d'},T')$.
\end{lemma}

\begin{proof}
Let $G=H_{d}\times\Z\times\S^1$, endowed with multiplication given by
$$
 (g,m,z)\cdot(g',m',z')=(gg'z'^m,m+m', zz')\ ,
$$
where we consider $z'^m$ as an element of $\S^1=(N_{d})_2$. 
Then $G$ is a 
$2$-step nilpotent Lie group, with 
$G_2=(H_{d})_2\times\{0\}\times\{1\}$.
The subgroup $\Gamma=\Lambda_{d}\times\Z\times\{1\}$ is a discrete 
cocompact subgroup of $G$. Write $X=G/\Gamma$,
$\sigma=(\tau,1,e(t))$, and let $S$ be the translation by $\sigma$ on 
$X$.

Note that $G$ is spanned by $G_0=H_{d}\times\{0\}\times\S^1$ and 
$\sigma$, and that the rotation induced by $\sigma$ on 
$G/\Gamma G_2=\T^{2d}\times \T$ is minimal by the hypothesis on $t$.
Therefore $(X,S)$ is 
minimal. It is an elementary nilsystem, written in reduced form.

The map $(g,m,z)\mapsto (g,z)$ from $G$ to 
$H_{d}\times\S^1$ induces a homeomorphism from $X$ onto 
$N_{d}\times\S^1$. We identify these spaces. 
The transformation $S$ has the form 
$(x,z)\mapsto (z\cdot Tx,ze(t))$. 
If $h$ is a function belonging to $\CC_1(N_{d},T)$, then the function
 $h\colon(x,z)\mapsto f(x)$ belongs to $\CC_1(X,S)$ and satisfies
$h(S^n(x,0))=q_n(t)f(T^nx)$ for every $x\in X$ and every $n\in\Z$. 
We deduce that $\Nil(X,S)$ is included in $\uq(t)\cdot\Nil(N_{d},T)$ 
and these spaces are equal by irreducibility 
(Corollary~\ref{cor:NXirreducible}). 
In particular, $\uq(t)\cdot\Nil(N_{d},T)$ does not contain the 
constant sequence $1$; 
substituting $-t$ for $t $ we have that $\uq(t)\notin\Nil(N_d,T)$.

We are left with showing that when $(N_{d'},T')$ is a minimal Heisenberg 
system, we have  $\Nil(N_{d'},T')\neq\uq(t)\cdot\Nil(N_d,T)$. 
Assume that these spaces are equal, that is, that 
$\Nil(X,S)=\Nil(N_{d'},T')$.
Let the subgroups $\Lambda'_1$ of $\Lambda_{d'}$  
and $\Gamma_1$ of $\Gamma$ 
 and the isomorphism $\Phi\colon \Lambda'_1\to\Gamma_1$ be given by 
Proposition~\ref{prop:equalortho2} and Lemma~\ref{lem:Phi}.
Since $\Gamma_1$ is of finite index in 
$\Gamma=\Lambda_{d}\times\Z\times\{1\}$, the restriction to $\Gamma_1$ 
of the natural projection of $\Gamma$ onto $\Lambda_d$ can not be one 
to one.  
Therefore, there exists $0\neq\gamma\in\Lambda'_1$ such that 
$\Phi(\gamma)=(0,m,1)$ for some non-zero $m\in\Z$. 
By property~\eqref{eq:Phicommut} of Lemma~\ref{lem:Phi}, we have that 
$$
 [\gamma,\tau']=[\Phi(\gamma),\sigma]= e(mt)\ .
$$
But $[\gamma,\tau']$ is an eigenvalue of $(N_{d'},T')$ and thus 
of $(N_d,T)$, a 
contradiction of the hypothesis.
\end{proof}

\begin{lemma}
\label{lemma:ouf}
Let $(N_d,T)$ and $(N_{d'},T')$ be two minimal Heisenberg systems of dimensions $2d+1$ and $2d'+1$ 
respectively, with translations given by
$$
 \tau=(\alpha_1,\dots,\alpha_d,\beta_1,\dots,\beta_d,e(\gamma))
\text{ and }
 \tau'=(\alpha'_1,\dots,\alpha'_{d'},\beta'_1,\dots,\beta'_{d'},
e(\gamma'))
$$
respectively. Let $t\in\T$.

Then $\uq(t)\Nil(N_d,T)=\Nil(N_{d'},T')$ if and only if $d=d'$ and there exist
a $2d\times 2d$ matrix $Q\in\Sp_{2d}(\Q)$ 
and integers $m,k_1,\dots,k_d,\ell_1,\dots,\ell_d$ with $m\geq 1$
such that 
\begin{align}
\label{eq:ouf}
 & Q\begin{pmatrix}
\alpha_1+k_1/m\\ \dots \\ \alpha_d+k_d/m\\ \beta_1+\ell_1/m
\\  \dots \\ \beta_d+\ell_d/m
\end{pmatrix}=\begin{pmatrix}
\alpha'_1\\ \dots \\\alpha'_d\\ \beta'_1\\ \dots \\ \beta'_d
\end{pmatrix}\\
\label{eq:ouf2}
\text{ and }\quad &
mt= \sum_{i=1}^d(k_i\beta_i-\ell_i\alpha_i)\bmod 1\ .
\end{align}
\end{lemma}

\begin{proof}\strut

\medskip\noindent a) 
Let $(N_d,T)$ be a minimal Heisenberg system of dimension $2d+1$, 
where $T$ is the translation by 
$\tau=(\alpha_1,\dots,\alpha_d,\beta_1,\dots,\beta_d,e(\gamma))$. 
Let $m,k_1,\dots,k_d,\ell_1,\dots,\ell_d$ be integers  with $m\geq 1$
and assume that $t\in\T$ satisfies~\eqref{eq:ouf2}.
Let 
$$
\tau'=(\alpha_1+k_1/m,\dots,\alpha_d+k_d/m,\beta_1+\ell_1/m,\dots,\beta_d+k_d/m,e(\gamma'))
$$
where $\gamma'\in\T$ is arbitrary and $T'$ the translation by $\tau'$ 
on $N_d$.
The first $2d$ coordinates of $\tau'$ are rationally 
independent modulo $1$ and thus $(N_d,T')$ is minimal. We 
compare $\Nil(N_d,T)$ and $\Nil(N_d,T')$.

We have that $\tau'^m=\tau^m u \gamma $, where $u$ is some element of 
$(H_d)_2=\S^1$ and 
$$
\gamma=(k_1,\dots,k_d,\ell_1,\dots,\ell_d,1)\in\Gamma\ .
$$
 An immediate computation gives that for every integer $p$,
$$
 \tau'^{mp}=\tau^{mp}v^{mp}e\bigl(\frac{mp(mp-1)}2t\bigr)\gamma^p
$$
for some $v\in (H_d)_2=\S^1$. 

Let $x_0$ be the image in $N_d$ of the unit element of $H_d$ and let $f$ be 
a function on $N_d$, belonging to $\CC_1(N_d)$, with $f(x_0)\neq 0$.
Let $\ua$ and $\ub$ be the  sequences given by  
$a_n=\one_{m\Z}(n)f(T^nx_0)$ and $b_n=\one_{m\Z}(n)f(T'^nx_0)$. 
 These sequence are
not identically zero, the sequence $\ua$ belongs to  $\Nil(N_d,T)$
and the sequence $\ub$ belongs to $\Nil(N_d,T')$ by definition. 
But for every integer $n$ we have $b_n=v^nq_n(t)a_n$ by the above 
computation and thus the sequence $\ub$ belongs to 
$\uq(t)\Nil(N_d,T)$.  Since the spaces $\uq(t)\Nil(N_d,T)$ and 
$\Nil(N_d,T')$ are irreducible (Corollary~\ref{cor:NXirreducible}), 
they are equal: $\uq(t)\Nil(N_d,T)=\Nil(N_d,T')$.

\medskip\noindent b)
Now assume that $d=d'$ and that $t\in\T$, 
the matrix $Q\in\Sp_{2d}(\Q)$,
and the integers $m,k_1,\dots,k_d,\ell_1,\dots,\ell_d$ satisfy the 
conditions~\eqref{eq:ouf} and~\eqref{eq:ouf2}. 
We need to show that 
$\uq(t)\cdot\Nil(H_d,T)=\Nil(H_d,T')$. 
By using Lemma~\ref{lem:class-Heisen}, we immediately reduce to the 
case that $Q$ is the identity matrix. 
The result then follows from part a) above.

\medskip\noindent c) Assume now that
$\uq(t)\cdot\Nil(H_d,T)=\Nil(H_d,T')$.  We show that there exist 
$m,k_1,\dots,k_d,\ell_1,\dots,\ell_d$ and $Q$ satisfying~\eqref{eq:ouf} 
and~\eqref{eq:ouf2}.

By Lemma~\ref{lem:convsnot}, there exists a non-zero integer $m$ such 
that $mt$ is an eigenvalue of $X$ and thus there 
exist integers $k_1,\dots,k_d,\ell_1,\dots,\ell_d$ such that 
relation~\eqref{eq:ouf2} holds. For $1\leq j\leq d$ we set 
$\alpha''_j=\alpha_j+k_j/m$ and $\beta''_j=\beta_j+\ell_j/m$. Let 
$T''$ be the translation by 
$\tau''=(\alpha''_1,\dots,\alpha''_d,\beta''_1,\dots,\beta''_d,\gamma)$ 
on $N_d$. By Part~a), $(N_d,T'')$ is minimal and 
$\Nil(N_d,T'')=\uq(t)\cdot\Nil(N_d,T)=\Nil(N_{d'},T')$. The 
conclusion follows from Lemma~\ref{lem:class-Heisen}.
\end{proof}

\subsection{Proof of Theorem~\ref{th:classification}}
\label{subsec:proof-classific}

\noindent a)
We check that if 
$t,d,\alpha_,\dots,\alpha_d,\beta_1,\dots,\beta_d$ are as in the 
statement of the theorem, then the sequence
$\ua=\uq(t)\uomega(\alpha_1,\beta_1).\dots.\uomega(\alpha_d,\beta_d)$ 
is an elementary nilsequence.

If $d=0$ and $t=0$, this sequence is constant. If $d=0$ and $t$ is 
irrational, then $\ua=\uq(t)$ and this sequence belongs to 
$\Nil(Y)$ for 
some affine elementary nilsystem (Section~\ref{subsec:affine}).

Assume now that $d>0$. 
By the discussion of Section~\ref{subsec:nilseqHeisen},
the sequence 
$\ub=\uomega(\alpha_1,\beta_1).\cdots.\uomega(\alpha_d,\beta_d)$ 
belongs to $\Nil(N_d,T)$ for some 
Heisenberg system $N_d$.
Let $t\in\T$ be an irrational  and let $(Y,T')$ be an  affine elementary 
nilsystem with $\uq(t)\in\Nil(Y,T')$. 
If for every integer $m\neq 0$, $mt$ does not belong to 
the eigenvalue group of $N_d$, then by Lemma~\ref{lem:convsnot} the 
sequence $\uq(-t)$ does not belong to $\Nil(N_d,T)$.  Thus 
$\Nil(Y,T')\neq\Nil(N_d,T)$. 
By Lemma~\ref{lem:convsnot}, there exists an
 elementary nilsystem $(X,S)$ with 
$\uq(t)\cdot\Nil(N_d,T)=\Nil(X,S)$. Therefore
 $\uq(t)\ub$ belongs to 
$\Nil(X,S)$ and is an elementary nilsequence. 

We are left with the case that $mt$ belongs to the group of 
eigenvalues of $(N_d,T)$ for some integer $m\neq 0$. Let 
$k_1,\dots,k_d,\ell_1,\dots,\ell_d$ be given by~\eqref{eq:ouf2}, $Q$ 
be the identity $2d\times 2d$ matrix  and 
set $\alpha'_j=\alpha_j+k_j/m$ and $\beta'_j=\beta_j+\ell_j/m$ for 
$1\leq j\leq d$. Then 
$(\alpha_1,\dots,\alpha_d,\beta_1,\dots,\beta_d)$ are rationally 
independent modulo $1$ and $N_d$ endowed with the translation$T'$  by 
$\tau'=(\alpha'_1,\dots,\alpha'_d,\beta'_1,\dots,\beta'_d)$ is a 
minimal Heisenberg system. By Lemma~\ref{lemma:ouf}, $\ua=\uq(t)\ub$ 
belongs to $\Nil(N_d,T')$ and thus is an elementary nilsequence.

\medskip\noindent b)
Let $\ua$ be a non-identically zero elementary nilsequence. We show 
that its class contains some sequence of the announced type.

If it is almost periodic, then its class is $\AP$ and thus contains the 
constant sequence $1$.  This is the exactly the sequence $\uq(0)$.
Now assume that $\ua$ is not almost periodic. 
By Lemma~\ref{lem:basicelement}, there exists a minimal nilsystem $X$ 
with $\CS(\ua)=\Nil(X)$. By Proposition~\ref{prop:non-connected}, 
there exist $t\in \T$ and an elementary nilsystem $(X'=G'/\Gamma',T')$ with 
$G'$ connected and $\Nil(X,T)=\uq(t)\cdot\Nil(X',T')$. If $t$ is  
rational, we have that $\uq(t)\cdot\Nil(X',T')=\Nil(X',T')$ 
and thus we can assume that $t$ is equal to $0$ or is irrational.
By Lemma~\ref{lem:connected-Heisen}, we can assume that $(X',T')$ is a
Heisenberg system and 
$\Nil(X,T)$ contains a sequence of the type 
$\uomega(\alpha_1,\beta_1).\cdots.\uomega(\alpha_d,\beta_d)$ 
satisfying the independence condition.
We conclude that $\CS(\ua)$ contains a 
sequence of the announced form.

\medskip\noindent c)
Let $t,t'\in\T$ be equal to $0$ or be irrationals, $d,d'\geq 0$ and let
 $\alpha_1,\dots,\alpha_d,\beta_1,\dots,\beta_d$ and
$\alpha'_1,\dots,\alpha'_{d'},\beta'_1,\dots,\beta'_{d'}$ satisfy
the independence condition.  Assume that the sequences
$\uq(t)\,\uomega(\alpha_1,\beta_1).\cdots.\uomega(\alpha_d,\beta_d)$ and 
$\uq(t')\,\uomega(\alpha'_1,\beta'_1).$ $\cdots.\uomega(\alpha'_{d'},\beta'_{d'})$ 
belong to the same class. Then the sequences
$\uq(t-t')\,\uomega(\alpha_1,\beta_1).\cdots$ $.\uomega(\alpha_d,\beta_d)$ and 
$\uomega(\alpha'_1,\beta'_1).$ $\cdots.\uomega(\alpha'_{d'},\beta'_{d'})$ belong to the same class.
 By the first part of 
Lemma~\ref{lem:convsnot}, $d$ and $d'$ are either both zero or both 
non-zero.  If $d=d'=0$, then $\uq(t-t')$ is almost periodic by 
Proposition~\ref{prop:equalortho2} and thus $t-t'$ is rational. If 
$d$ and $d'$ are non-zero, then Lemma~\ref{lemma:ouf} gives the 
announced relations between the parameters 
$t,t',d,d',\alpha_1,\dots,\beta'_d$.

The converse implication also follows from the same lemmata.
\qed

\subsection{Proof of Theorem~\ref{th:span}}
The theorem follows immediately from Theorem~\ref{th:classification} 
and Proposition~\ref{prop:decomposition}.\qed

\subsection{Proof of Theorem~\ref{th:span2}}
\label{subsec:proofspan2}
\begin{lemma}
\label{lem:quadratic}
Let $\ua$ be a sequence belonging to the family $\CM$. 
Then the linear span of $\{\ue(u)\ua\colon u\in\T\}$ is dense in 
$\CS(\ua)$ for the quadratic norm.
\end{lemma}

\begin{proof}[Proof of Lemma~\ref{lem:quadratic}]
If $\ua=\ue(s)\uq(t)$ for some $s, t\in\T$, then one can easily check 
that $\CS(\ua)=\AP\cdot\ua$.  Moreover, the linear span of 
 $\{\ue(u)\ua\colon u\in\T\}$ is dense in 
$\CS(\ua)$ for the uniform  norm.  Thus it is dense for the 
quadratic norm.

Now consider the case that 
$\ua=\ue(s)\uq(t)\ub$ where 
$\ub=\uomega(\alpha_1,\beta_1).\cdots.\uomega(\alpha_d,\beta_d)$ 
with $d\geq 1$ and the independence condition 
satisfied.   Since
$\CS(\ua)=\uq(t)\CS(\ub)$, we can reduce to the case that 
$s=t=0$, meaning that $\ua=\ub$.

In Section~\ref{subsec:nilseqHeisen}, we 
built  an explicit  minimal Heisenberg system 
$(N_d,T)$, $x_0\in N_d$  and an explicit function $f\in\CC_1(N_d)$ 
 such that $b_n=f(T^nx_0)$ for every $n\in\Z$.

Let $K=\{w\in N_d\colon f(w)=0\}$. Since $f\in\CC_1(N_d)$, 
$K=\pi\inv(L)$ where $\pi$ is the natural projection from $N_d$ onto 
$Z=H_d/(H_d)_2\Lambda_d$ and $L$ is a closed subset of $Z$.
Recall that $Z$ can be identified with $\T^{2d}$.

The relation~\eqref{eq:omegatheta} between the function 
$\kappa$ defined in~\eqref{eq:omega} and the theta function shows that 
$\tilde f(x,y,z)$ defined by~\eqref{eq:fomega} vanishes only when 
one of the coordinates $x_j$ or $y_j$ of $x$ or $y$ is equal to 
$1/2$ modulo $1$. Therefore, the subset $L$ of $Z=\T^{2d}$ has zero 
Haar measure. It follows 
(as in the proof of Corollary~\ref{cor:NXirreducible}) 
that the linear span of the family of functions
$\{f\cdot\chi\circ\pi\colon \chi\in\wh Z\}$ is dense in $\CC_1(N_d)$ for 
the $L^2(N_d)$ norm. 

For $\chi\in\wh Z$ and for every $n\in\Z$, we have 
$$
(f\cdot\chi\circ\pi)(T^nx_0)=e(ns)f(T^n x_0)=e(ns)b_n \ ,
$$ 
where $s\in\T$ is 
the eigenvalue associated to $\chi$. Let $E$ be the group of 
eigenvalues of $T$. We deduce that the closure in quadratic norm of the
linear span of the family
$\{\ue(s)\ub\colon s\in E\}$ of sequences contains all sequences of 
the form $(h(T^nx_0)\colon n\in \Z\bigr)$ for $h\in\CC_1(N_d)$.
Therefore, the closure in quadratic norm of the family
$\{\ue(s)\ub\colon s\in \T\}$ of sequences contains 
all sequences of 
the form $(e(nt)h(T^nx_0)\colon n\in \Z\bigr)$ for $h\in\CC_1(N_d)$ 
and $t\in\T$.
By Lemma~\ref{lem:basepoint}, it therefore contains $\CS(\ub)$.
The announced result 
follows.
\end{proof}

\begin{proof}[Proof of Theorem~\ref{th:span2}]
The result follows immediately from 
Proposition~\ref{prop:decomposition},
Theorem~\ref{th:classification} and Lemma~\ref{lem:quadratic}.
\end{proof}

\section{Characterization of elementary nilsequences}
\label{sec:charactelemen2}
In this Section we show Theorem~\ref{th:charactelemen}.
\subsection{The ``only if'' part}
\label{sec:charactelemen1}
First we show the easy part of
Theorem~\ref{th:charactelemen}.  We restate it here for 
convenience.
\begin{proposition*}
Let $\ua$ be an elementary nilsequence. 
There exists a compact (in the norm topology) subset $K\subset 
\ell^\infty(\Z)$ such that for all $k\in \Z$, there exists $t\in\T$ 
such that the sequence $\ue(t)\sigma^k\ua=\bigl(e(nt)a_{n+k}\colon 
n\in\Z\bigr)$  belongs to $K$.
\end{proposition*}

\begin{proof}
If $\ua$ is almost periodic, then the result follows immediately from 
the characterization~\eqref{item:third} of almost periodic sequences 
in Proposition~\ref{prop:AP}.  Thus we assume that $\ua$ is not 
almost periodic.

By Lemma~\ref{lem:basicelement}, there exists an elementary nilsystem 
$(X,T)$ with $\ua\in\Nil(X)$. Let $(\tilde X,\tilde T)$ be the 
Bohr extension of $(X,T)$ as defined in 
Section~\ref{subsec:constr-universal} and let $\tilde x_0\in\tilde X$.
By Lemma~\ref{lem:nilXtilde}, there exists a continuous function
$f$ on $\tilde X$ with $a_n=f(\tilde T^n\tilde x_0)$ for all $n\in\Z$.

Recall that $\tilde X$ can be written as $\tilde G/\tilde\Gamma$, 
where $\tilde G$ is a locally compact $2$-step nilpotent group,
$\tilde\Gamma$ is a  closed cocompact subgroup, and $\tilde T$ is 
translation by some $\tilde\tau\in\tilde G$.
The natural projection $\tilde G\to\tilde X$ is an open map and by 
compactness, we deduce that there 
exists a  compact subset $H$ of $G$ such that $G = H\Gamma$. 
Let $K$ be the subset of $\ell^\infty(\Z)$ 
consisting of sequences of the form 
$$
\bigl( f(g\tilde \tau^n \cdot\tilde x_0)\colon n\in \Z
\bigr)
$$ 
for $g\in H$.  Since the map $(\tilde g,\tilde x)\mapsto \tilde 
g\cdot\tilde x$ is continuous, the function $(g,x)\mapsto f(g\cdot 
x)$ is uniformly continuous on $H\times X$ and 
 $K$ is a compact subset of $\ell^\infty(\Z)$.

Let $k$ be an integer.  Pick $g\in H$ and $\gamma\in \Gamma$
with $\tau^k = g\gamma$.  The sequence 
$\ub = \bigl(f(g\tau^n\cdot x_0) \colon n\in\Z\bigr)$ belongs to $K$.  
Taking $t\in\T$ with  $e(t) = [\gamma, \tau]\in G_2=\S^1$, we have that 
for all $n\in\Z$, 
$$
(\sigma^k\ua)_n = a_{n+k} = f(\tau^k\tau^n\cdot x_0) 
= f(g\gamma\tau^n\cdot x_0) = e(nt)f(g\tau^n\cdot c_0) \ .
$$
Thus $\sigma^k\ua = \ue(t)\ub$ and this completes the proof.  
\end{proof}

\subsection{The ``if'' part}
In the next  Subsections we show the implication
\eqref{it:compact} $\Longrightarrow$ \eqref{it:elementary}
of Theorem~\ref{th:charactelemen}. We restate it for convenience.

\begin{proposition*}
Let $\ua=(a_n\colon n\in\Z)$ be a bounded sequence and assume 
that there exists a compact 
(in the norm topology) subset $K$ of
$\ell^\infty(\Z)$ such that for all $k\in \Z$, there exists $t\in T$ 
such that the sequence $\ue(t)\sigma^k\ua=\bigl(e(nt)a_{n+k}\colon 
n\in\Z\bigr)$  belongs to $K$. 
Then $\ua$ is an elementary 
nilsequence.
\end{proposition*}
We use the following notation. Recall that $\norm\cdot_\infty$ 
denotes the uniform norm on $\ell^\infty(\Z)$.
 Let $B$ be a subset of $\ell^\infty(\Z)$.  
For $t\in\T$,  $\ue(t)\cdot B:=\{ \ue(t)\ub\colon \ub\in B\}$.
If $k\in\Z$, $\sigma^k B:=\{\sigma^k\ub\colon \ub\in B\}$.
 If $z$ is a number, $z B=\{z\ub\colon \ub\in B\}$.

Before the proof, we need a simple Lemma:

\begin{lemma}
\label{lem:syndetic}
Let $E$ be  a syndetic subset of $\Z$ of syndetic constant $L$, meaning that
$$
 \bigcup_{j=0}^{L-1}(j+E)=\Z\ .
$$
If $t\in\T$ is such that 
$$
|e(nt)-1|<1/3L\text{ for every $n\in E$ and all }|t|<1/3L\ ,
$$
then $t=0$.
\end{lemma}
\begin{proof}[Proof of the lemma]
If $|e(nt)-1|<1/3L$ for all $n\in E$, then for every $n\in\Z$ 
there exists $j$ with $0\leq j<L$ such that $n\in E+j$.  Thus
the point $nt$ belongs to the interval $(jt-1/3L,jt+1/3L)$.

Therefore, $\Z t$ is included in the union $U$ of the intervals 
$(jt-1/3L,jt+1/3L)$ for $0\leq j<L$. 
The complement of $U$ in $\T$ 
contains a closed interval $J$ of length
$1/3L$.  As the sequence $(nt\colon n\in\Z)$ avoids $J$,
$|t|\geq 1/3L$ and we have a contradiction.
\end{proof}

\subsection{First reductions}

We now turn to the proof of the Proposition.
We can obviously assume that the sequence $\ua$ is not almost periodic 
and in particular not identically zero. We can also assume that
$$\norm\ua_\infty=1\ .$$
Assume that $\ua$ belongs to $K$.
Substituting  the compact set
$$
\S^1\cdot K:=\{u\ub\colon u\in\S^1,\ \ub\in K\}
$$
for $K$, we can assume that $K$ is invariant under multiplication by 
constants of modulus $1$. Write
$$
 \Omega=\{u\ue(t)\cdot \sigma^n\ua\colon u\in\S^1,\ t\in\T, \ n\in\Z\}
$$
and let $\bar\Omega$ denote the norm closure of $\Omega$ in 
$\ell^\infty(\Z)$. We remark that $\bar\Omega$ is invariant under the 
shift $\sigma$, under multiplication by constants of modulus $1$ 
and under the operators of multiplication by exponential sequences.
By hypothesis, we have
$$
 \Omega\subset\bigcup_{t\in\T} \ue(t)\cdot K\ .
$$
Substituting $K\cap\bar\Omega$ for $K$ we can assume that
$$
 K\subset\bar\Omega\ .
$$

\subsection{First step} We obviously have that 
$$
 \norm{\ub}_\infty=1\text{ for every }\ub\in\bar \Omega\ .
$$

\begin{claim}
\label{cl:one}
$\bar\Omega$ is closed under pointwise convergence.
\end{claim}

\begin{proof}
Let $(\ub_j\colon j\in\N)$ be a sequence in $\Omega$, converging 
pointwise to a sequence $\ub$. 
For every $j\in\N$, we write
$$
 \ub_j=\ue(t_j)\uc_j
\text{ where }t_j\in\T\text{ and }\uc_j\in K\text{ for every }j\ .
$$
By passing to a subsequence, we can assume that the sequence 
$(t_j\colon j\in\N)$ converges to some $t$ in $\T$ and that 
 the sequence $(\uc_j\colon j\in\N)$ converges uniformly to some $\uc\in K$.  Thus
$\ub=\ue(t)\uc$.
Since 
$\bar\Omega$ is invariant under multiplication by $\ue(t)$ and contains $K$,
we have that $\ub\in\bar\Omega$.
\end{proof}

\begin{claim} 
\label{cl:two}
The set $S=\{n\in\Z\colon |a_n|\geq 2/3\}$ is syndetic.
\end{claim}

\begin{proof}
Assume that $S$ is not syndetic. 
Then there exists a sequence $(k_j\colon j\in\N)$
 of integers such that $(\sigma^{k_j}\ua)$ converges pointwise to a 
sequence $\ub$ with $\norm{\ub}_\infty\leq 2/3$. 
By Claim~\ref{cl:one}, $\ub\in\bar\Omega$, and thus $\norm{\ub}_\infty=1$, 
a contradiction.
\end{proof}

In the sequel, we write $L$ for the syndetic constant of the set 
$S$.
  
From the last claim we immediately 
deduce:
\begin{claim} 
\label{cl:twobis}
For every $\ub\in\bar\Omega$, the set $\{n\in\Z\colon |b_n|\geq 1/2\}$
is syndetic of syndetic constant $L$.
\end{claim}

\begin{claim}
\label{cl:three}
For every compact subset $Q$ of $\bar\Omega$, 
the family of sets $\bigl\{\ue(t)\cdot Q\colon t\in\T\}$ is locally finite.
More precisely, every ball $B$ of radius $r<1/25\,L$ 
intersects
$\ue(t)\cdot Q$ for  only finitely many values of $t$.
\end{claim}

\begin{proof}
Let $B$ be a  ball of radius $r<1/25\,L$ and let $(t_j\colon j\in\N)$ 
be a sequence 
of distinct elements of $\T$ such that  $B\cap \ue(t_j)\cdot Q\neq\emptyset$
for every $j\in\N$. 

For every $j$, we chose $\uc_j\in Q$ with $\ue(t_j)\uc_j\in B$.
 Passing to a subsequence, we can assume that the 
sequence $(t_j\colon j\in\N)$ 
converges to some $t$ in $\T$ and that the sequence 
$(\uc_j\colon j\in\N)$ converges uniformly to some sequence $\uc\in Q$. 
We note that $\uc\in \bar\Omega$.

For all $i,j\in\N$, by the triangle inequality,
\begin{multline*}
\norm{\ue(t_j-t_i)\uc}_\infty
\leq \norm{\ue(t_j)\uc-\ue(t_j)\uc_j}_\infty+ 
\norm{\ue(t_j)\uc_j-\ue(t_i)\uc_i}_\infty +
\norm{\ue(t_i)\uc-\ue(t_i)\uc_i}_\infty\\
\leq\norm{\uc-\uc_j}_\infty+\norm{\uc-\uc_i}_\infty+2r
\end{multline*}
because  $\ue(t_i)\uc_i$ and $\ue(t_j)\uc_j$ belong 
to $B$.

Let $i$ and $j$  be sufficiently large.
The last quantity is $<4r$ and thus 
$$
|e((t_j-t_i)n)-1|< 8r<1/3L
$$ 
for $n$ in the set
$F=\{n\in\Z\colon|c_n|\geq 1/2\}$. As we also have that 
$|t_i-t_j|<1/3L$, Claim~\ref{cl:twobis} and Lemma~\ref{lem:syndetic} 
give that $t_i=t_j$, a contradiction.
\end{proof}

We immediately deduce:

\begin{claim}
\label{cl:four1}
For every compact subset $Q$ of $\bar\Omega$, the union of the sets $\ue(t)\cdot Q$ 
for $t\in \T$ is closed.
\end{claim}

\begin{claim}
\label{cl:four}
$$
 \bar\Omega=\bigcup_{t\in\T}\ue(t)\cdot K\ .
$$
\end{claim}
\begin{proof}
The union on the right hand side is closed by Claim~\ref{cl:four1}, 
it contains $\Omega$ and is included in $\bar\Omega$ by construction.
\end{proof}

\begin{claim}
\label{cl:five}
$\bar\Omega$ is locally compact.
\end{claim}
\begin{proof}
Every closed ball of radius $r<1/25L$ is covered by sets of the form 
$\ue(t)K$ by Claim~\ref{cl:four} and intersects only finitely many of them by
Claim~\ref{cl:three}. It is thus  included in a finite union of  sets 
of this form and so is compact.
\end{proof}

\subsection{A group of transformations}
Let $\Isom(\bar\Omega)$ denote the group of isometries of 
$\bar\Omega$,  endowed with the topology of pointwise 
convergence.  This topology coincides with the compact-open topology, 
as well as with the topology of pointwise convergence on the dense  subset 
$\Omega$ of $\bar\Omega$.  We write the
action of this group on $\Omega$ as 
$(g,\ub)\mapsto g\cdot\ub$ for $g\in\Isom(\bar\Omega)$ and 
$\ub\in\bar\Omega$.
Let $G$ be the closure in $\Isom(\bar \Omega)$ of the group spanned 
by the shift and the operators of multiplication by $\ue(t)$ for $t\in 
\T$. This operator is  written simply $\ue(t)$.

We remark that for $t\in\T$, the commutator of $\sigma$ and $\ue(t)$ is 
multiplication by the constant $e(-t)$.
Therefore,  the commutator subgroup of $G$ is 
the group of multiplication by constants of modulus $1$, identified 
with $\S^1$, and $G$ is $2$-step nilpotent.

\begin{claim}
\label{cl:eight}
Let $(u_i\colon i\in I)$ be a generalized sequence (a filter) in $\S^1$,
$(t_i\colon i\in I)$ a generalized  sequence in $\T$, and
$(k_i\colon i\in I)$  a generalized  sequence in $\Z$.  
For every $i\in I$, let 
$g_i= u_i\ue(t_i)\sigma^{k_i}\in G$. Then the generalized  sequence 
$(g_i\colon i \in I)$ 
converges in $G$ if and only if
\begin{enumerate}
\item
The generalized  sequence $(g_i\cdot\ua\colon i\in I)$ 
converges in $\bar\Omega$;
\item
\label{it:tit}
The generalized  sequence $(u_i\colon i\in I)$ converges in $\S^1$;
\item\label{it:dual}
The generalized sequence $(k_i\colon i\in I)$ 
converges in the Bohr group $\BZ$, 
meaning that for 
every $t\in D$ the generalized  sequence $(k_it\colon i\in I)$ 
converges in $\T$.
\end{enumerate}
\end{claim}

\begin{proof}
Let $\ub=u\ue(t)\sigma^k\ua\in\Omega$, with $u\in\CS^1$, $t\in\T$, 
$k\in\Z$. For every $i\in I$ we have 
$$
 g_i\cdot\ub=e(kt_i-k_it)\,u\ue(t)\sigma^k\, g_i\cdot\ua\ .
$$
The three conditions of the claim are equivalent to the property that 
$g_i\cdot\ub$ converges for every $\ub\in\Omega$. 
\end{proof}
\begin{claim}
\label{cl:nine}
$G$ acts transitively on $\bar\Omega$.
\end{claim}
\begin{proof}
Let $\ub\in\bar\Omega$. 
Since $\Omega$ is dense in $\bar\Omega$, there exist sequences 
$(s_i\colon i\in I)$, $(t_i\colon i\in I)$ and $(k_i\colon i\in I)$ 
as in Claim~\ref{cl:eight} such that 
that $e(s_i)\ue(t_i)\sigma^{k_i}\cdot\ua\to\ub$. Substituting 
subsequences for these sequences we can assume that  
properties~\eqref{it:tit} and~\eqref{it:dual} are satisfied. 
The limit $g$ of the sequence $(e(s_i)\ue(t_i)\sigma^{k_i})$ satisfies 
$g\cdot\ua=\ub$.
\end{proof}

\begin{claim}
\label{cl:ten}
For every compact subset $L$ of $\bar\Omega$, the subset 
$\tilde L=\{g\in G\colon g\cdot\ua\in L\}$ of $G$ is compact.
Moreover, $G$ is locally compact.
\end{claim}
\begin{proof}
Let $\ub\in\Omega$. Choose  $h\in G$ such that 
$h\cdot\ua=\ub$. 
For every $g\in\tilde L$, we have $g\cdot\ub=ghg\inv h\inv hg\cdot a$.  
Thus $\tilde L\cdot \ub$ is included in the 
compact set $\S^1h\cdot L$ and so is relatively compact. The first 
part of the Claim follows from  Ascoli's Theorem, and the second part 
from the fact that $\bar \Omega$ is locally compact 
(Claim~\ref{cl:five}).
\end{proof}

\subsection{The space $X$}
Let  $\Gamma$ be the subgroup $\{\ue(t)\colon t\in \T\}$ of 
$G$.
It follows immediately from
 Claim~\ref{cl:three} that $\Gamma$ is 
discrete and closed in $G$.

Write $X=G/\Gamma$ endowed with the quotient topology.
For $h\in G$, the map $g\mapsto hg$ from $G$ to $G$ induces a homeomorphism 
of $X$ that we write $x\mapsto h\cdot x$. 
We remark that the map $G\times X\to X$ defined by $(h,x)\mapsto 
h\cdot x$ is continuous.

By Claim~\ref{cl:ten}, $\tilde K=\{g\in G\colon g\cdot\ua\in K\}$ is 
a compact subset of $G$ and by Claim~\ref{cl:four} we have that
$G$ is the union $\Gamma \tilde K$ of the sets 
$\ue(t)\cdot\tilde K$ for 
$t\in  \T$. Since $K$ is invariant under multiplication by constants of 
modulus one, $\tilde K$ is invariant under multiplication by $\S^1$ 
and $\Gamma\tilde K=\tilde K\Gamma$. Therefore the restriction to 
$\tilde K$ of the natural projection $G\to X=G/\Gamma$ is onto and we 
have:

\begin{claim}
\label{cl:eleven}
$X$ is compact.
\end{claim}

Let the function $\tilde f$ on $G$ be defined by 
$$
\tilde f(g)=\bigl(g\inv\cdot\ua\bigr)_0\ ,
$$
where the subscript $0$ denotes the coordinate corresponding 
to $0$.
This function is continuous and for every $t\in\T$,
$$
\tilde 
f(g\ue(t))=\bigl(\ue(-t)g\inv\cdot\ua\bigr)_0=
\bigl(g\inv\cdot\ua\bigr)_0=\tilde f(g)\ .
$$
(Recall that $\ue(-t)\ub$ is the product of the sequence $\ub$ and the 
sequence $\ue(t)=(e(nt)\colon n\in\Z)$.  Therefore we have
$(\ue(-t)\ub)_0=b_0$.)
Thus the function $\tilde f$ induces a continuous function $f$ 
on $X$. We remark that for every $x\in X$ and every $u\in\S^1$, we have
$f(u\cdot x)=uf(x)$. 
Using terminology similar to that in Definition~\ref{def:elementary} 
and Section~\ref{subsec:uses-for-elementary}, we say that 
the function $f$ belongs to  $\CC_{1}(X)$.

Let $T\colon X\to X$ be the map $x\mapsto \sigma\inv\cdot x$ and let $x_0$ be 
the image in $X$ of the unit element of $G$.
For every $n\in\Z$, we have
$$
 f(T^n x_0)=\tilde f(\sigma^{-n})=
\bigl(\sigma^n\cdot\ua)_0=a_n\ .
$$

\subsection{Conclusion}
This situation was studied in [HM] in the case that $G$ is metrizable, 
 but the same proof extends to the present case. We outline the method.

It is classical~\cite{MZ} that the locally compact $2$-step nilpotent 
group $G$ is an inverse limit of $2$-step nilpotent Lie groups. This
means that there exists a decreasing family $(K_i)$ of compacts 
subgroups of the center of $G$, with trivial intersection, such that 
$G/K_i$ is a $2$-step nilpotent Lie group for every $i$. Therefore 
the system $(X,T)$ is the inverse limit of the $2$-step nilsystems
$(X_i=G/K_i\Gamma, T_i)$ where $T_i$ is the transformation of $X_i$ 
induced by $T$. For every $i$, the commutator subgroup 
of $G/K_i$ is the circle group, and thus $(X_i,T_i)$ is an elementary 
nilsystem. For every $i$, let $p_i\colon X\to X_i$ be the natural 
projection and set $x_i=p_i(x_0)$. The function $f$ on $X$ is the uniform 
limit of functions of the form $f_i\circ p_i$, where 
$f_i\in\CC_1(X_i)$ for every $i$. This gives us that 
the sequence $\ua$ is the 
uniform limit of the sequences  $\ua(i):=(f_i(T_i^nx_i)\colon 
n\in\Z)$, and so is an elementary nilsequence.\qed

\appendix

\section{Proof of Lemma~\ref{lem:tildeX}}
\label{appendix:A}

We recall the notation introduced in Section~\ref{subsec:constr-universal}.
We have that $(X=G/\Gamma,T)$  is an elementary nilsystem, where $T$ is the 
translation by $\tau\in G$ and $\mu$ is its Haar measure. The
Kronecker factor of this system is $(Z,m_Z,S)$ and $\pi\colon X\to Z$ 
is the factor map. As before, $p\colon G\to Z$ is the natural homomorphism, so 
that $S$ is the translation by $\sigma=p(\tau)$ on $Z$. 
Finally, $\BZ$ is the Bohr compactification of $\Z$, endowed with the 
translation $R$ by $1$ and with its Haar measure $m$, and  $r\colon\BZ\to Z$ 
is the continuous homomorphism characterized by $r(1)=\sigma$.

The system $(\tilde X,\tilde T)$ denotes the Bohr extension
of $(X,T)$: $\tilde X=\{(x,z)\in X\times\BZ\colon\pi(x)=r(z)\}$ 
endowed with the translation $\tilde T=T\times R$; $q_1\colon\tilde 
X\to X$ and $q_2\colon\tilde X\to\BZ$ are the natural factor maps.

We first show that the system $(\tilde X,T)$ is 
uniquely ergodic and that the invariant measure is the conditional 
independent joining (we recall the definition below).

\subsection{Notation}
Let $\nu$ be a $\tilde T$-invariant measure on $\tilde X$. 

By unique ergodicity of $(\BZ, R)$ and $(X,T)$, the images of 
$\nu$ under $q_1$ and $q_2$  are equal to 
 $m$ and $\mu$, respectively  (In other words, $\nu$ is a 
joining). 
We write $\E_r\colon L^2(m)\to L^2(m_Z)$ and $\E_\pi\colon 
L^2(\mu)\to L^2(m_Z)$ for the operators of conditional expectation, 
defined as usual by:
\begin{multline*}
\text{for every }\phi\in L^2(m),\ \psi\in L^2(\mu)\text{ and }
h\in L^2(m_Z),\\
\int_Z\E_r \phi\cdot h\,dm_Z=\int\phi\cdot h\circ r\,dm
\text{ and }
\int\E_\pi\psi\cdot h\,dm_Z=\int \psi\cdot h\circ p\,d\mu\ .
\end{multline*}
We claim that $\nu$ is the conditional independent joining 
$\mu\times_Z m$ of $\mu$ 
and $m$ over $Z$, meaning it is the measure on $X\times \BZ$ defined by:
\begin{multline*}
\text{for every }\phi\in L^2(\mu) \text{ and }\psi\in L^2(m),\\
\int\psi(x)\psi(w)d \mu\times_Z m(x,w)
=\int \E_p\phi(z)\cdot\E_r\psi(z)\,dm_Z(z)\ .
\end{multline*}

\subsection{Proof of the Claim}
The measure $\nu$ defines a continuous operator 
$\Phi\colon L^2(m)\to L^2(\mu)$ by:
\begin{multline*}
 \text{for every }
\phi\in L^2(\mu)\text{ and every }\psi\in 
L^2(m),\\
\int \phi(x)\psi(w)\,d\nu(x,w)=\int \Phi\psi(x)\, 
\phi(x)\,d\mu(x)\ .
\end{multline*}
Since $\nu$ is invariant under $\tilde T=T\times R$, we have that 
$\Phi\circ R=T\circ \Phi$. 
If $\psi$ is a character of $\BZ$, it is an eigenfunction of 
$(\BZ,R)$.  Thus  $\Phi\phi$  is an 
eigenfunction of $X$ and by hypothesis $\Phi\psi$ factorizes through 
$Z$.
Since the characters of $\BZ$  span $L^2(m)$, the same 
property holds for every $\psi\in L^2(m)$. 

Therefore the last integral is equal to
$$
 \int \Phi\psi(x)\,\E_p\phi\circ p(x)\,d\mu(x)=
\int \psi(w)\,\E_p\phi\circ p(x)\,
d\nu(x,w)\ .
$$
since $\nu$ is concentrated on $\tilde X$, it follows that $p(x)=r(w)$ 
for $\nu$-almost every $(x,w)$ and the last integral can be 
rewritten as
$$
\int \psi(w)\,\E_p\phi\circ r(w)\,
d\nu(x,w)=
\int \psi(w)\,\E_p\phi\circ r(w)\, dm(w)\ .
$$
The projection of $\nu$ on $Z$ is equal to $m$ and using the
definition of the operator $\E_p$, this is equal to
$$
 \int \E_p\phi(z)\,\E_r\psi(z)\,dm_Z(z)= 
\int \phi(x)\,\psi(w)\,d\mu\times_Z m(x,w) \ .
$$
By definition of this measure, and our claim is proved.

\subsection{} In particular, $(\tilde X,\tilde T)$ is uniquely 
ergodic with invariant measure $\mu\times_Z m$. 
The topological support 
of  this measure 
 is clearly equal to $\tilde X$.
\qed

\section{Elementary nilsystems arising from connected groups}
\label{sec:desc-conected}
We explain how elementary nilsystems with a 
connected group have a simple form.  
Similar results appear in~\cite{Lei2}.  
Simply connected, connected nilpotent Lie groups admitting a discrete 
cocompact subgroup have been well understood since 
Malcev~\cite{M} and we refer to Chapter 5 
of~\cite{CG} for their properties. 
The proof reduces to this case and we only give a sketch of the reduction. 

\begin{lemma}
\label{lem:desc-connected}
Let $(X=G/\Gamma,T)$ be an elementary nilsystem (minimal, written in 
reduced form) with connected group $G$. 
Then $(X,T)$ is the product of a rotation on a finite dimensional torus 
with an elementary nilsystem $(X'=G'/\Gamma',T')$ of the following type:

$d\geq 1$ is an integer, $A$ is a $2d\times 2d$ matrix with integer 
entries such that the matrix $B=A-A^t$ is nonsingular; 
$G=\R^{2d}\times\S^1$ with multiplication given by
$$
(x,z)\cdot(x',z')=\bigl(x+x',zz'e(\langle Ax\mid x'\rangle)\bigr)
$$
and $\Gamma=\Z^d\times\{1\}$.
\end{lemma}

\begin{proof}
Let $\tilde G$ be the simply connected covering of $G$, $p\colon 
\tilde G\to G$ be the natural projection and
let  $\tilde\Gamma=p\inv(\Gamma)$. 
Then $\tilde G$ is a simply connected, 
connected $2$-step nilpotent Lie group, the kernel of $p$ is a 
discrete subgroup of the center of $\tilde G$, and $\tilde\Gamma$ is a 
discrete cocompact subgroup of $\tilde G$.
Since $G_2$ is the circle group, $\tilde G_2$ is of dimension $1$. 

Let $W_1,\dots,W_m,Z$ be a Malcev base of the second kind of the Lie algebra of
$\tilde G$. 
This means that every $g\in \tilde G$ can be written in a unique 
way as 
$$
g=\exp(w_1W_1)\dots\exp(w_mW_m)\exp(zZ)\text{ with }
w_1,\dots,w_m,z\in\R\ ,
$$
that $\tilde G_2=\{\exp(zZ)\colon s\in\R\}$,
and that 
$$
\tilde\Gamma=\{ 
\exp(p_1X_1)\dots\exp(p_mX_m)\exp(qW)\colon p_1,\dots,p_m,q\in\Z\}\ .
$$


We use the
coordinates of $G$ associated to this base
and identify $\tilde G$ with $\R^m\times\R$, $\tilde 
G_2$ with $\{0\}\times\R$ and $\tilde\Gamma$ with $\Z^m\times\Z$.
Then the multiplication in $\tilde G$ is given by
$$
 (w_1,\dots,w_m,z)\cdot(w'_1,\dots,w'_m,z')=
(w_1+w'_1,\dots,w_m+w'_m,z+z'+\langle Aw\mid w'\rangle)
$$
where $w=(w_1,\dots,w_m)$, $w'=(w'_1,\dots,w'_m)$ and  $A$ is an $m\times m$ 
matrix with integer entries. The commutator map is given by 
$$
 \bigl[(w_1,\dots,w_m,z)\cdot(w'_1,\dots,w'_m,z')\bigr]=
(0,\dots,0 ,\langle Bw\mid w'\rangle)\ ,
$$
where $B=A^t-A$. 


We can write $m=2d+q$ and find a base of $\Z^{2d}\times\Z^q$ of the form 
$(e_1,\dots,e_{2d},f_1,\dots,f_q)$ such that $(f_1,\dots,f_q)$ is a 
base of the kernel of the matrix $B$. The center of $\tilde G$ is thus 
the linear span of the vectors $f_1,\dots,f_q$ and 
$\exp(Z)=(0,\dots,0,1)$.

By hypothesis, the intersection of the center of $G$ and $\Gamma$ is 
trivial. It follows that $f_1,\dots,f_q$ and $\exp(Z)$ belong to the 
kernel of the projection $p$.

Let $\Lambda$ be the subgroup spanned by these elements, 
$G' =G/\Lambda$ and let $\Gamma'$ be the image of $\tilde\Gamma$ in 
$G'$. Then $X=G'/\Gamma'$ and $G',\Gamma'$ have the announced form.
\end{proof}

\section{Non-examples of nilsequences}
\label{ap:non-example}

\subsection{A non-almost periodic example} 
First we show that the sequence $\ub$ 
defined by $b_n=e(\lfloor n\alpha\rfloor \beta\colon n\geq 1)$, 
where $\alpha, \beta$ are independent over the rationals, 
is not almost periodic.  

Let $(k_j\colon j\geq 1)$ 
be a sequence of integers such that the sequence $(k_j\alpha\colon 
j\geq 1)$ converges
to $0$ modulo $1$ and $(\lfloor k_j\alpha\rfloor\beta\colon j\geq 1)$ 
converges to some $\gamma$ modulo $1$.

If $\ub$ is almost periodic, then by passing to a subsequence, 
which we also denote by $(k_j\colon j\geq 1)$, 
we can assume that the sequence $(\sigma^{k_j}\ub\colon j\geq 1)$ 
converges uniformly to some sequence $\uc$. 
On the other hand, it is easy to check that $(\sigma^{k_j}\ub\colon j\geq 1)$ 
converges to $e(\gamma)\ub$ in quadratic norm, but not uniformly.

We have that $\uc=e(\gamma)\ub$, a contradiction.

\subsection{A $2$-step nonexample}
We show that the sequence $e(\lfloor n\beta\rfloor n\alpha)$, 
where $\alpha, \beta$ are independent over the rationals, 
 is not a 
$2$-step nilsequence.  

We proceed by contradiction.  Assume that 
this sequence is a $2$-step nilsequence. 
Then using Proposition~\ref{prop:closed}, the sequence $\ua$ given by
$$
 a_n=e(\lfloor n\beta\rfloor n\alpha-\frac {n(n-1)}2\alpha\beta )
$$
also is a $2$-step nilsequence.

Consider the Heisenberg nilmanifold $N_1 = G/\Gamma$, 
endowed with translation by $\tau=(\alpha,\beta,0)$. 
Let $e$ denote the image in $N_1$ of the unit element of $G$. 
The subset 
$F=[0,1)\times[0,1)\times\T$ of $H_1=\R\times\R\times\T$ is a 
fundamental domain for the projection $p\colon H_1\to N_1$.
Thus there exists a unique function $h$ on $N_1$ such that 
$h(p(x,y,z)=e(-z)$ for every $(x,y,z)\in F$. 
Note that $h$ is Riemann integrable, but is not 
continuous, and that for every $n$ we have
$a_n=h(T^ne)$.
(In fact a general statement along these lines holds. 
Let $(X,T)$ be a minimal $k$-step nilsystem, $x_0\in X$, and $f$ be
a Riemann integrable function on $X$.  If the 
sequence $(f(T^nx_0)\colon n\geq 1)$ is a $k$-step nilsequence,
then $f$ is continuous.)

Let $(h_j\colon j\geq 1)$ be a sequence of continuous functions on 
$N_1$ converging to $h$ in $L^2(N_1)$ and for every $j$, let
$\ua(j)$ denote the sequence $(h_j(T^ne)\colon n\in\Z)$. 
By Proposition~\ref{prop:elementary-closed} and 
Lemma~\ref{lem:basepoint}, each of 
these sequences is an elementary nilsequence belonging to 
$\Nil(N_1)$.  
As $j\to+\infty$, $\ua(j)$ converges to $\ua$ in quadratic norm 
(see Section~\ref{subsec:first}).

By Proposition~\ref{prop:decomposition} 
and Theorem~\ref{th:equalortho}, proceeding as in the proof of 
Lemma~\ref{lem:basicelement}, 
it is easy to deduce that if a nilsequence is a limit in 
the quadratic norm of elementary nilsequences belonging to $\Nil(X)$, then 
it also belongs to $\Nil(X)$.  In particular, it is an elementary 
nilsequence.  Therefore 
$\ua$ is an elementary nilsequence.

Let $(k_j\colon j\geq 1)$ be a sequence of integers such that 
$k_j\alpha$ and $k_j\beta$ tend to $0$ modulo $1$ and 
$\frac{k_j(k_j-1)}2\alpha\beta$ converges modulo $1$ to some $\gamma$.

By Theorem~\ref{th:charactelemen}, 
passing to a subsequence of $(k_j\colon j\geq 1)$ if necessary, we 
can assume that there exists a sequence $(t_j\colon j\geq 1)$ 
in $\T$ such that the sequence
$(\ue(t_j)\sigma^{k_j}\ua\colon j\geq 1)$ 
converges uniformly as $j\to\infty$ to some sequence.

On the other hand, it is easy to check that 
that the sequence 
$(\ue(k_j\alpha\beta)\sigma^k\ua\colon j\geq 1)$ 
converges as $j\to\infty$ in quadratic norm, but not 
uniformly. We deduce that the sequence $(t_j-k_j\alpha\beta\colon 
j\geq 1)$ is eventually constant, and 
a contradiction follows.


\begin{thebibliography}{99999}

\bibitem{AGH}
L.~Auslander, L.~Green and F.~Hahn.
Flows on homogeneous spaces.
{\em Ann. Math. Studies}
{\bf 53}, Princeton Univ. Press (1963).


\bibitem{BHK}
V.~Bergelson, B.~Host and B.~Kra, with an Appendix
by I.~Ruzsa.  Multiple recurrence
and nilsequences.  {\em Inventiones Math.} {\bf 160} (2005), 261-303.

\bibitem{BL}
V.~Bergelson and A.~Leibman.  
Distribution of values of bounded generalized polynomials. 
{\em Acta Math.} {\bf 198} (2007) 155-230.
\bibitem{CG}
L.~Corwin and F.P.~Greenleaf,
\emph{Representations of nilpotent Lie groups and their applications},
Cambridge University Press (1990).

\bibitem{F}
H.~Furstenberg.  
Ergodic behavior of diagonal measures and a theorem of Szemer\'edi
  on arithmetic progressions.
{\em J. d'Analyse Math.} {\bf 31}  (1977), 204--256.  

\bibitem{GT} B.~Green and T.~Tao.  Linear equations in the primes.
To appear, {\em Annals of Math}.

\bibitem{GT2} B.~Green and T.~Tao.  Quadratic uniformity of the 
M\"obius function.  To appear, {\em Annales de l'Institut Fourier}.

\bibitem{HK6}
B.~Host and B.~Kra.
Uniformity seminorms on $\ell^\infty(\Z)$ and an inverse 
theorem.  In preparation. 


\bibitem{HM} B.~Host and A.~Maass.
Nilsyst\`emes d'ordre deux et parall\'el\'epip\`edes.
To appear, {\em Bull. Math. Soc. France}.

\bibitem{Lei2}
A.~Leibman.
Pointwise convergence of ergodic averages for polynomial sequences of
rotations of a nilmanifold.  
{\em Ergod. Th. \& Dynam. Sys.} {\bf 25} (2005), 201-113.

\bibitem{Les}
E.~Lesigne.
Sur une nil-vari\'et\'e, les parties minimales
associ\'ees \`a une translation sont uniquement ergodiques.
{\em Ergod. Th. \& Dynam. Sys.}
{\bf 11} (1991), 379--391.

\bibitem{M} A.~Malcev.  On a class of homogeneous spaces. 
{\em Amer. Math. Soc. Transl.} {\bf 9} (1962), 276-307.

\bibitem{MZ} D.~Montgomery and L.~Zippin.  {\em Topological Transformation 
Groups}. Interscience Publishers (1955). 

\bibitem{Pa1}
W.~Parry.
Ergodic properties of affine transformations and flows on
nilmanifolds.
{\em Amer. J. Math.}
{\em 91} (1969), 757--771.

\bibitem{Pa2}
W.~Parry.
Dynamical systems on nilmanifolds.
{\em Bull. London Math. Soc.}
{\bf 2} (1970), 37--40.


\bibitem{R}
D.~J.~Rudolph.
Eigenfunctions of $T\times S$ and the
Conze-Lesigne algebra.
{\em  Ergodic Theory and its Connections
with Harmonic Analysis}, Eds. K.~Petersen and I.~Salama.
Cambridge University Press,
New York (1995), 369--432.


\end{thebibliography}
\end{document}